\documentclass[12pt]{article}
\usepackage{philhspstyle}

\begin{document}

\title{The entropy efficiency of point-push mapping classes on the
punctured disk}

\author{Philip Boyland and Jason Harrington}

{\namefont

\maketitle

}




\begin{abstract}
We study the maximal entropy per unit generator of push-point
mapping classes on the punctured disk. Our work is motivated
by  fluid mixing by rods in a planar domain.
If a single rod moves among $N$-fixed obstacles, the resulting
fluid diffeomorphism is in the push-point mapping class
associated with the loop in $\pi_1(D^2 - \{N\ \text{points}\})$
traversed by the single stirrer. The collection of motions in each
of which the stirrer goes 
around a single obstacle generate the group of push-point
mapping classes, and the entropy efficiency with respect to
these generators gives a topological measure of the 
mixing per unit energy expenditure of the mapping class.
We give lower and upper bounds for $\Eff(N)$, the maximal
efficiency in the presence of $N$ obstacles, and prove
that  $\Eff(N)\raw \log(3)$ as $N\raw\infty$.
For the lower bound we compute the entropy efficiency
of a specific push-point protocol, $\HSP_N$, which we
conjecture achieves the maximum. The entropy computation
uses the action on chains in a
$\Z$-covering space of the punctured disk 
which is designed for push-point protocols.
For the upper bound we estimate the exponential growth rate
of the action of the  push-point mapping classes on the fundamental 
group of the
punctured disk using a collection of incidence matrices 
and then computing the generalized spectral radius of the
collection.  
\end{abstract}

\section{Introduction}
The topological entropy $h(f)$ of a map
$f$  is a standard measure of its dynamical
complexity. The entropy of an isotopy class is 
the infinimum of all the entropies of the \homeos\ in the class.
If a surface isotopy class has positive entropy,
a theorem of Thurston guarantees the existence of
an (almost) canonical map in the isotopy class which
realizes the infinimum and whose dynamics 
are described by a finite collection of mixing subshifts
of finite type (\cite{thurston}, \cite{FLP}). 
Further, Handel's 
theorem shows that these dynamics are present up to
semiconjugacy in every \homeo\ in the mapping class
 (\cite{handel}, \cite{bdisostab}).

The prototypical positive entropy mapping class is
of \pA\ type and the Thurston representative in
the class is a \pA\ map  $\phi$. In this case,
$\lambda = \exp(h(\phi))$ is the dilation of the
\pA\ map, an algebraic integer with a wide variety of geometric
and topological interpretations. As a consequence
of these interpretations 
it is natural to seek the \pA\ map with least dilation 
within semi-groups of isotopy classes (see, for example,
\cite{minent1},  \cite{minent2}, or \cite{minent3}).
From a dynamics perspective, maximizing complexity is
also of interest. Since entropy is multiplicative
under iteration, $h(f^n) = n h(f)$, the entropy over a
semigroup  is always unbounded, and thus one needs a
normalization.
 In many cases the semi-group of interest has a natural set of generators,
for example, a  collection of Dehn twists
or motions of stirring rods as discussed below.
In these cases following  Thiffeault/Finn and Moussafir
 (\cite{finnth2}, \cite{otherone}, \cite{finnth1}) 
we consider the \textit{entropy efficiency}
of a mapping class which is the topological entropy per
unit generator. The entropy efficiency of the semigroup
is then the supremum of the efficiencies of its members.

As a consequence of Thurston's theorem, the
entropy of an isotopy class is equal to the exponential
growth rate of words under iteration by the induced map
on the fundamental group of the surface. One may also
define this growth rate and thus the entropy
of any automorphism of a finitely generated group $G$. 
Given a finite collection of automorphisms $\cA = \{a_1, \dots, a_n\}$
of the group,
we define the \de{entropy efficiency} of the semigroup
generated by $\cA$ as the maximal entropy per unit generator
(see \S\ref{enteffsec}). If $G$ is a free Abelian group, the automorphisms
are linear isomorphisms and  the entropy efficiency of
$\cA$ is the log of the generalized spectral radius of
the collection of linear transformations (Remark~\ref{gsree}).
The generalized and joint spectral radii are 
very useful  and much-studied linear objects, and the entropy efficiency
is their non-linear analog.

Our main motivation for studying the entropy efficiency
is fluid mixing via the motion of stirring rods in a planar
domain (\cite{BAS},   \cite{finnth1}, \cite{finnth2}, \cite{boylandeuler}). 
In the simplest  situation the mixing is accomplished
by  $n$-rods which after time-$T$  
return to their original positions. Such a stirring protocol
is naturally described as a braid with $n$-strings. 
When the rods move through the fluid there is a resulting
motion of the entire fluid region, and the evolution of a fluid particle
from its initial position to its position after time-$T$
 defines the \de{\fmh}.
The exact nature of this \diffeo\ depends on the physical properties
of the fluid and the geometry of the fluid region and
rod paths. However, as long as the fluid
 does not tear, all possible \fmh s arising from the
same stirring protocol are in the
same mapping class, and this class is the one usually
associated with the braid (see \S\ref{braidsec}). Thus if we can
determine the entropy of this isotopy class we have
a lower bound for the complexity of the evolution of
any fluid motion resulting from the stirring protocol.
The topological entropy has a particular physical
significance  because it measures, among other things, 
the exponential growth rate of material
lines in the fluid under the action of the protocol.
Exponential stretching of material lines does
not guarantee good mixing, but in many situations 
this rate of stretching is a good topological
``first order'' indicator.

The other consideration in implementing 
a stirring protocol is its energy cost.
Protocols which maximize the mixing per unit expenditure
of energy could have great practical value.
he energy required to implement a stirring protocol certainly
depends on the geometry of the path of the stirrer, but since we are
characterizing  protocols just by their braids,
 we seek a topological measure of effort which just
depends on the braid.
Thus we pick some class of fundamental motions of the stirrers,
for example, the switching of a pair of rods, and treat
these as generators of a semigroup of mapping
classes. The corresponding maximal entropy efficiency
 gives a rough notion of the maximal mixing
per unit expenditure of energy possible in the class
of protocols.  A protocol which realizes the maximum
would be optimal in a precise mathematical
sense, but perhaps not in any true physical sense, see
\cite{finnth1}.

The entropy efficiency of the full braid group with
respect to the standard generators was studied
by Finn and Thiffeault in \cite{finnth1}. 
They show that the maximal efficiency
 is  bounded above by $\log[(1 + \sqrt{5})/2]$. Further, 
when $n > 4$,  this bound is strict, and    
when $n = 3$ and  $4$ the bound is realized by $\sigma_1 \sigma_2\I$
and $\sigma_1\I \sigma_2 \sigma_3\I \sigma_2$, respectively. 
They also studied a generating set consisting
of words of the form $\sigma_{i_1}^{\pm 1} \dots
 \sigma_{i_k}^{\pm 1}$, where for all $j$,
$|i_j - i_{j+1}| > 1$. For this class they
obtained $\log(1 + \sqrt{2})$ as an upper bound for the
maximal efficiency. 
This class of generators models  motions
of the collection of rods in which the entire motion
can be done simultaneously. 

Here we study the class of stirring protocols in which there
are $N$-fixed obstacles and a single stirrer moves.
Thus the stirring protocols are determined by the 
path of the single stirrer in the $N$-punctured disk.
This path can be considered an element of the fundamental
group and so these stirring motions are called 
\de{$\pi_1$-stirring protocols}.
Each $\pi_1$-stirring protocol determines  a braid on $(N+1)$ strings and
an isotopy class on the $(N+1)$-punctured disk.
The mapping classes which result from $\pi_1$-stirring protocols
are called \de{point-push} in the mathematics literature
(\cite{birmanbook}, \cite{kra}, \cite{farbstudent}, \cite{primer}).
Thus the  collection of $\pi_1$-stirring protocols, 
denoted here as  $\cP_N$,
is the subgroup of the braid group
on $(N+1)$ strings corresponding to the point-push
isotopy classes in the $N$-punctured disk. To maintain the
traditional mathematical terminology while incorporating
the physical  intuition and motivation   of our problem, we will usually
use  \de{\posp} to refer to   elements of  $\cP_N$.

 We use the analog of
the standard generators of the fundamental group of the
$N$-punctured disk as our fundamental stirring motions.
Thus  for $j = 1, \dots, N$, let $\alpha_j$ correspond
to moving the  stirrer counterclockwise around the 
$j^{th}$-obstacle (Figure~\ref{alphafig}),  and 
we consider the entropy efficiency of a \posp\ 
with respect to the  generating set $\{\alpha_1, 
\dots, \alpha_N\}$ of $\cP_N$. The maximal
entropy efficiency over $\cP_N$ is denoted    
 $\Eff(N)$.
Our main theorem is:
\begin{theorem}\label{main}
 For $N\geq 2$ there are explicitly defined matrices $\HN$ and
$\hHN$ such that 
the maximal efficiency $\Eff(N)$ of a \posp\  
in the presence of $N$ obstacles satisfies
\begin{equation*}
\frac{\log(3^N - 3N - 1) - \log(N)}{N}
 \leq \frac{\log(\rho(\HN))}{N} \leq \Eff(N) \leq 
\frac{\log(\rho(\hHN))}{N} \leq \frac{\log(3^N -2)}{N},
\end{equation*}
and so $\Eff(N) \raw \log(3)$ as $N\raw\infty$.
\end{theorem}
Note that the lower bound  has no content when $N=2$.
The most efficient protocol with two obstacles is given
in Theorem~\ref{casetwothm}.

`
The proof of Theorem~\ref{main}
  follows directly from the lower bound given
in Theorem~\ref{lowerboundthm} and 
the upper bound given in Theorem~\ref{upperboundthm}.
Both the upper and lower bound proceed by reducing
the problem to the computation of 
the spectral radius of a specific product of generating
matrices, but the reduction in the two cases is quite different.
In both cases the required computations are tractable
because of the simple forms of actions made possible by
the adapted set of generators shown in Figure~\ref{genfig}.

Since the maximal efficiency is the supremum over the
efficiencies of all \posp s, the efficiency of any
single protocol serves as a lower bound.
In \S\ref{hspdefinesec} we define a \posp\ $\HSP_N$ 
which numerically appears
to give the maximal efficiency (Figure~\ref{hspfig}). The 
invariant train track for $\HSP_N$ is rather simple,
but the corresponding Markov transition matrix has very
large entries. Thus the computation or even estimation 
of the dilation of $\HSP_N$ as its  spectral radius
appears to be quite difficult. Instead we use homological
methods developed by Fried (\cite{friedtwisted}, \cite{friedzeta},
\cite{bandboyland}).
 In short, the
idea is to lift the mapping class to finite covers where
the dilation is computable as the spectral radius of the
action on homology. This cover is found as the quotient
of a $\Z$-covering space whose homology is treated as
a free $\Z[t\pmo]$-module.  The action of the mapping
class is then a matrix over $\Z[t\pmo]$, and substitutions
of $t = \exp(2 \pi i p/q)$ give pieces of the action
of the mapping class on the homology of the corresponding
$\Z_q$-covering space.

As noted above, the entropy of a \posp\
is equal to the exponential growth rate of words under
iteration by the induced map
on the fundamental group of the punctured disk, $\pi_1(D_{N+1})$.
For each generator $\alpha_j$ of $\cP_N$, we find  the
incidence matrix $A^{(j)}$ of its action on $\pi_1(D_{N+1})$.
The incidence matrix just counts the    
number of occurrence of a generator in the image of another
generator. As also remarked above, the entropy efficiency is
the nonlinear analog of the \gsr\ and in 
Proposition~\ref{gsrupper} we show that for a finitely
generated semi-group of free group automorphisms, 
the log of the
\gsr\ of the  incidence matrices of the generators 
is an upper bound for the entropy efficiency of the semigroup. 
In general, finding a \gsr\ is difficult, but because of
the particularly simple form of the action of the $\alpha_j$
on the generators in Figure~\ref{genfig}, in our
case we can determine the exact product which yields
the \gsr\ and estimate its spectral radius (Proposition~\ref{gsrcompute}).

We conjecture that
the protocol $\HSP_N$ which is used for the lower
bound actually realizes the maximal efficiency of a \posp\
in the presence of $N$ obstacles when $N\geq 3$.
 Figure~\ref{hspfig}(b) shows
a symmetric view of the path of the single stirrer 
of $\HSP_4$. The path lies on a hypotrochoid which
lead to the terminology  hypotrochoid stirring protocol.
The fact that the path of the single stirrer can be
generated by a disk rolling inside   a circle  
 has practical implications;  it implies that
the stirring protocol  may be implemented
using a fairly simple system of gears to move the stirrer,
see \cite{Kobayashi}, \cite{finnth1}. It is also interesting to note that
the spectral radii of the matrices used in Theorem~\ref{main}
appear to be algebraic units of a special type. Numerical computations
for $N\leq 20$ indicate that each $\rho(\HN)$ is a
Salem number and  each $\rho(\hHN)$ is a Pisot number. These
classes of algebraic numbers arise in a variety of topological
and dynamical circumstances, especially
in optimization problems (see, for example,
\cite{hironaka} or \cite{minent3}). The generator matrices 
whose product defines $\hHN$ are the absolute values of those
which define  $\hHN$; the connection of this fact to the algebraic
type of their spectral radii warrants further investigation.

\medskip
\noindent \textbf{Acknowledgement:} We would like to
thank Jean-Luc Thiffeault for suggesting this problem
and many useful conversations and thank Toby Hall for his assistance
with train track calculations.

\section{Preliminaries}
\subsection{Linear Algebra and the generalized spectral radius}
We recall some basic facts and establish a few conventions. 
The identity matrix of the appropriate dimension in
a given context is
denoted $I$. An inequality  of a pair of vectors
means that the inequality
holds in each component; the analogous convention holds for
inequalities of matrices of the same dimensions. 

For a square matrix $M$, $\rho(M)$ will denote its spectral radius and 
for $1\leq p   \leq \infty$, $\|M\|_p$ is  the matrix norm induced by
the vector norm $\|\vv\|_p$.
Recall the basic facts that $\rho(M) \leq \|M\|_p$ and
 $\|M\|_2 = (\rho(M M^T))^{1/2}$. Gelfand's formula says
that $\rho(M) = \lim_{n\raw\infty} \|M^n\|_p^{1/n}$.
This implies that  if $M$ and 
$M'$ are non-negative, square matrices with $M \geq M'$, then
$\rho(M) \geq \rho(M')$.

For a given matrix $M$, we let $c(M)$ be the row vector of
its column sums, and so $c_j(M) = \sum_i M_{ij}$. An easy
computation yields that for any pair of square matrices,
$c(M M') = c(M) M'$. Also note that for a positive matrix
$M$, $\max_j(c_j(M)) = \| M \|_1$.

We shall also work with free modules over the ring of
 all Laurent polynomials
in the variable $t$ with integer coefficients, $\Z[t^{\pm 1}]$.
 The collection of
 $(n\times n)$-matrices over $\Z[t^{\pm 1}]$
is denoted $\Mat(n, \Z[t^{\pm 1}])$.

We recall the definition of the \gsr.  For more information
see \cite{gsrbook}.
For a  finite set of matrices $\cM =  \{M_1, \dots, M_n\}$ and a
$k>0$ let 
\begin{equation}\label{gsrdef}
\rho_k(\cM) = \max\{\rho(M_{i_1} M_{i_2} \dots M_{i_k})^{1/k} : 
M_{i_j}\in \cM \}.
\end{equation}
Since norms are submultiplicative,
\begin{equation}\label{maxeq}
\rho_k(\cM) \leq \max\{\|M_{i}\|_p: 1 \leq i \leq n\}.
\end{equation} 
Thus the \de{generalized spectral radius} defined as
\begin{equation*}
\rho(\cM) = \limsup_{k\raw\infty} \rho_k(\cM)
\end{equation*}
is always finite. 
Note that since the spectral radius is
multiplicative, for all $n>0$, $\rho_{nk}(\cM) \geq \rho_k(\cM)$,
and thus 
\begin{equation}\label{supok}
\rho(\cM) = \sup\{ \rho_k(\cM): k\in\Z^+\}.
\end{equation}

\subsection{Word growth and the efficiency 
of automorphisms}\label{enteffsec} 
Let $G$ be a finitely generated group with 
generators $x_1, \dots, x_m$.  For $g\in G$, let 
$\ell(g)$ be the minimum word length in the generators.
If $a:G\raw G$ is an automorphism, its 
\de{entropy}  is defined as 
\begin{equation}\label{expgrowth}
h(a) = \max_i \limsup_{n\raw\infty} 
\frac{\log(\ell(a^n(x_i)))}{n}.
\end{equation}
Note that $h(a)$ does not depend on
the choice of generators of $G$.  
If $\cA\subset Aut(G)$ is a semigroup generated by 
$\{a_1, \dots, a_k\}$, for $a\in \cA$ define its 
\de{efficiency} with respect to the given generating set as 
\begin{equation*}
\eff(a;\{a_1, \dots, a_k\}) = \frac{h(a)}{\Ell(a)},
\end{equation*}
where $\Ell(a)$ is the minimal word length of $a$
 in the generators $\{a_i\}$.
 Finally, the \de{maximal efficiency} 
of $\cA$ with respect to the given generating set is 
\begin{equation*}
\Eff(\cA, \{a_1, \dots, a_k\}) = 
\sup\{\eff(a;\{a_1, \dots, a_k\}) : a\in \cA\}. 
\end{equation*}
Note that $\eff(a)$ and $\Eff(\cA)$ usually \textit{do}
depend on the choice of generators of $\cA$. 

\begin{remark}\label{gsree}
 If the group $G$ is free Abelian, then
an automorphism $A$ is a linear transformation and
treating $g\in G$ as a vector $\vg$, its word length 
is $\ell(g) = \|\vg \|_1$.  
It then follows  that $h(A) = \log(\rho(A))$. 
Further, as a consequence of \eqref{supok}, the efficiency
of the semigroup $\cA$ generated by a collection of linear 
isomorphisms $\{A_1, \dots, A_n\}$ is the log of their 
generalized spectral radius, 
\begin{equation*}
\Eff(\cA; \{A_1, \dots, A_n\} = 
\log(\rho(\{A_1, \dots, A_n\})).
\end{equation*}
Thus the entropy efficiency may be viewed as the nonlinear analog of
the generalized spectral radius.
\end{remark}

\subsection{Braids}\label{braidsec}
There are three of the many interpretations  of braids 
which will be used here: as the motion
of a collection of points in the disk, as  physical
braids, and as elements of a mapping class group.
 We describe each of
these briefly and informally. For more information see 
\cite{birmanbook}, \cite{primer}.
Fix a model of the two-dimensional disk $D^2$ and
$n$ distinguished points $\{p_1, \dots, p_n\} \subset Int(D^2)$.
The disk punctured at the distinguished points is
$D_n := D^2 - \{p_1, \dots, p_n\}$.

A motion of the distinguished points is described by
a family of paths $\eta_i:[0,1] \raw D^2$ for $i =1, \dots, n$
with $\eta_i(t) \not=\eta_j(t)$ for all $t\in [0,1]$ and
$i\not= j$. To find the corresponding physical braid 
  define arcs $\heta:[0, 1] \raw D^2 \times [0,1]$
by  $\heta_i(t) = (\eta_i(t), 1-t)$. The physical
braid corresponding to the motion of points is the
image of $(\heta_1, \dots, \heta_n)$. We can define
the product of two point motions by following one by the
other. This corresponds to stacking the 
associated physical braids one below the other. 

The elements of Artin's braid group on $n$-strings, denoted $B_n$,
correspond to equivalence classes of these objects. It
can be considered as homotopy classes of point paths or else
of physical braids. There are $(n-1)$-standard generators of
$B_n$ with $\sigma_i$ corresponding to switching
the $i^{th}$ and  $(i+1)^{th}$ strings (or points)
 counterclockwise.  
There are two types of relations which define the braid group:
$\sigma_i \sigma_{i+1} \sigma_i = 
\sigma_{i+1} \sigma_{i} \sigma_{i+1}$
for all $i$, and $\sigma_i \sigma_{j} = \sigma_j \sigma_i$ for 
$|i-j| > 1$.

For the connection to a mapping class group, we let
$\MCG(D_n)$ be the collection of isotopy classes of \homeos\
of $D_n$, where the \homeos\ and isotopies
are required to fix the boundary of the disk pointwise.
The collection of  mapping classes $\MCG(D_n)$ becomes a group
under composition. For each generator $\sigma_i$ of $B_n$,
let $f_i$ be a \homeo\ that switches $p_i$ and $p_{i+1}$
counterclockwise. The map $\sigma_i \mapsto [f_i] \in \MCG(D_n)$
then yields an isomorphism $B_n \isomorphic \MCG(D_n)$.

The \de{pure braid group} on $n$ strings is the subgroup $P_n$
of $B_n$ corresponding to paths $\eta_j$ with $\eta_j(0) =
\eta_j(1)$. Thus  the corresponding physical
 braids do not permute their endpoints and the mapping classes
in $D_n$ fix all the punctures. 

\medskip

\noindent\textbf{Conventions:} There are many possible conventions
for the braid group and its various interpretations. Since there
is not much standardization, we summarize our conventions.  
Braid words are written from left to right, and so $\sigma_1 \sigma_2$
means perform $\sigma_1$ and then $\sigma_2$. In the corresponding
physical braid, we draw $\sigma_2$ below $\sigma_1$. We maintain
the same conventions when treating braids as mapping classes,
and so $\sigma_1 \sigma_2$ means switch the first two punctures
from the left counterclockwise and then switch the second
 two punctures. Thus
the convention is the opposite of what one would use
when composing functions. To be consistent with these conventions,
when we have a matrix representation of the braid group, the
 matrices will act on a row vector from the right. Thus
the matrices are the transpose of the most common convention.
Specifically, the components of the image of the 
$i^{th}$-basis vector gives the $i^{th}$-row of the matrix.

\subsection{Thurston-Nielsen theory and the 
entropy of a mapping class}\label{thurstonsec}
Thurston's theorem classifies mapping classes on a surface $M$
into three classes: \pA, finite order and reducible (see 
\cite{thurston}, \cite{FLP}, \cite{casson}, \cite{primer}).
 In each case the class contains what is now called the
Thurston representative which has the same name, so
\pA\ classes contain \pA\ \homeos, etc.
 We define the (topological)
entropy of a surface mapping class $\chi$ as 
\begin{equation}\label{classent}
h(\chi)  = \inf\{\htop(f) : f\in\chi\}.
\end{equation}
Part of Thurston's classification theorem says that 
this infimum is achieved by the
Thurston representative in the mapping class.
In addition, the exponential
word growth under iteration by the mapping class
on $\pi_1(M)$ is equal to $h(\chi)$.  Specifically,
if $M$ is a surface with chosen basepoint $Q\in M$,
pick  a \homeo\
$f\in\chi$ with $f(Q) = Q$ and let $f_\sharp$ denote
its action on $\pi_1(M, Q)$. Then treating $f_\sharp$
as an automorphism of $\pi_1(M, Q)$ and defining 
$h(f_\sharp)$ as in \eqref{expgrowth}, Thurston's theorem says that 
$h(f_\sharp) = h(\chi)$. For a  \pA\ class $\chi$, the number 
$\lambda = \exp(h(\chi)) >1$ is called the \de{dilation}

We call a braid $\beta$ \pA, finite order, or reducible
according to whether its corresponding mapping
class has that type, and we define its entropy $h(\beta)$ 
as in \eqref{classent}, again using its
corresponding mapping class.

\section{Point-push  protocols and their topological efficiency}
In accordance with the motivation from fluid mechanics
we will sometimes use the terminology ``stirring protocol''
to refer to a braid or equivalently, a mapping class
on the punctured disk. We study a special class of 
stirring protocols here in which just one stirrer moves.
The other $N$ stirrers (or punctures) are stationary and will be referred
to as the \de{obstacles}. 
The initial and final point of the path of the single stirrer
is denoted $P$, and so the protocol is given by a  
 path $\Gamma:[0,1]\raw D_N$ with $\Gamma(0) = 
\Gamma(1) = P$. We may extend $\Gamma$ to an isotopy
$f_t$ on $D_N$ so that $f_t(P) = \Gamma(t)$ for all
$t\in [0,1]$. We further puncture $D_N$ at 
$P$ forming $D_{N+1}$, and    class of $f_1$ in 
$\MCG(D_{N+1})$ is the point-push class corresponding
to $\Gamma$.  Since $\Gamma$ is a closed loop starting
and ending at $P$ it gives an element
of the fundamental group with basepoint $P$, $\pi_1(D_N, P)$.
A theorem of Kra \cite{kra} and
Birman \cite{birmanbook}  says that the point-push mapping class 
is uniquely determined by this 
element in the fundamental group
and so  elements of this class of stirring motions are 
sometimes called \de{$\pi_1$-protocols}.

\subsection{The subgroup of \posp s}
Each \posp\ with $N$ obstacles yields a mapping class that fixes the
punctures and so they correspond to elements of the pure braid
group on $(N+1)$ strings. The collection of such classes is
a subgroup we denote $\cP_N$. If $h:P_{N+1} \raw P_N$ is the
homomorphism which deletes the left-most string of the pure
braid, then
$\cP_N = \ker(h)$. Equivalently, the \posp s  correspond 
to mapping classes on $D_{N+1}$ with the property that
when the left-most  puncture is filled in, the resulting
class in $\MCG(D_N)$ is the identity.

\begin{figure}[ht]
\begin{center}
\includegraphics[width=.30\linewidth,angle=0]{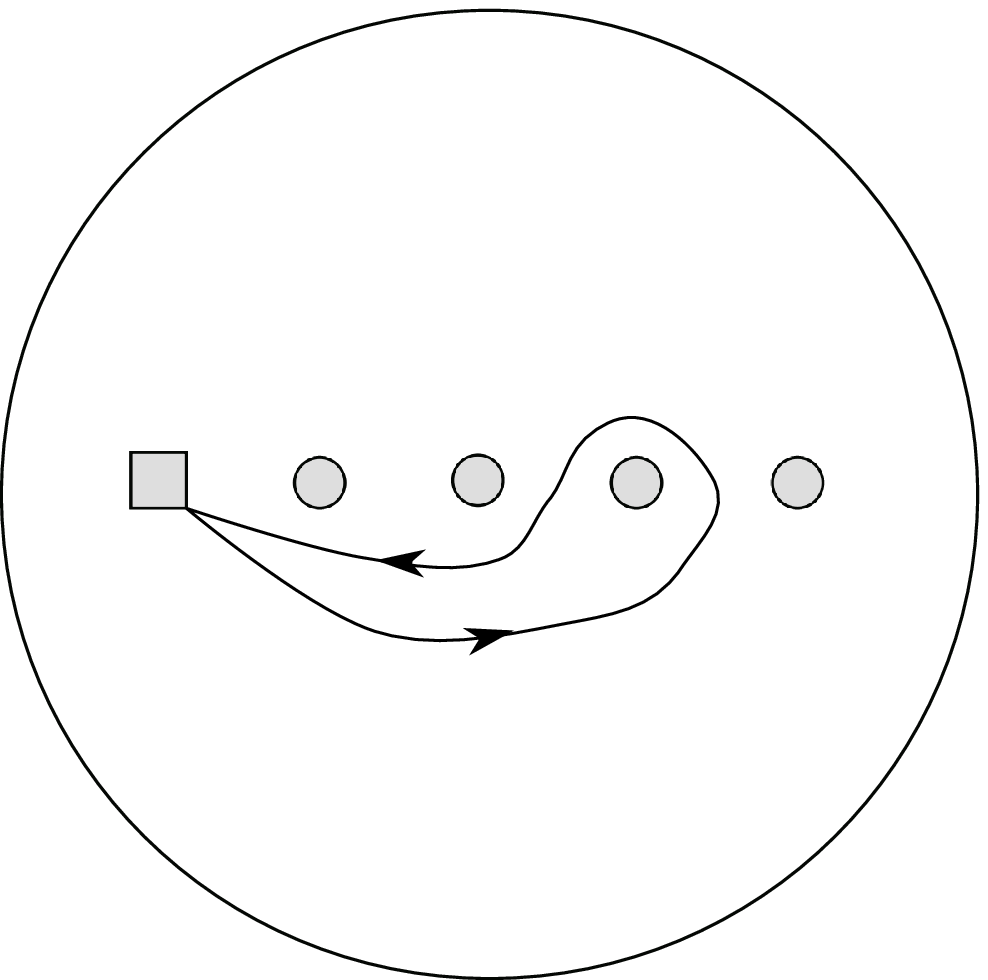}
\hskip 0.75in
\includegraphics[width=.18\linewidth,angle=0]{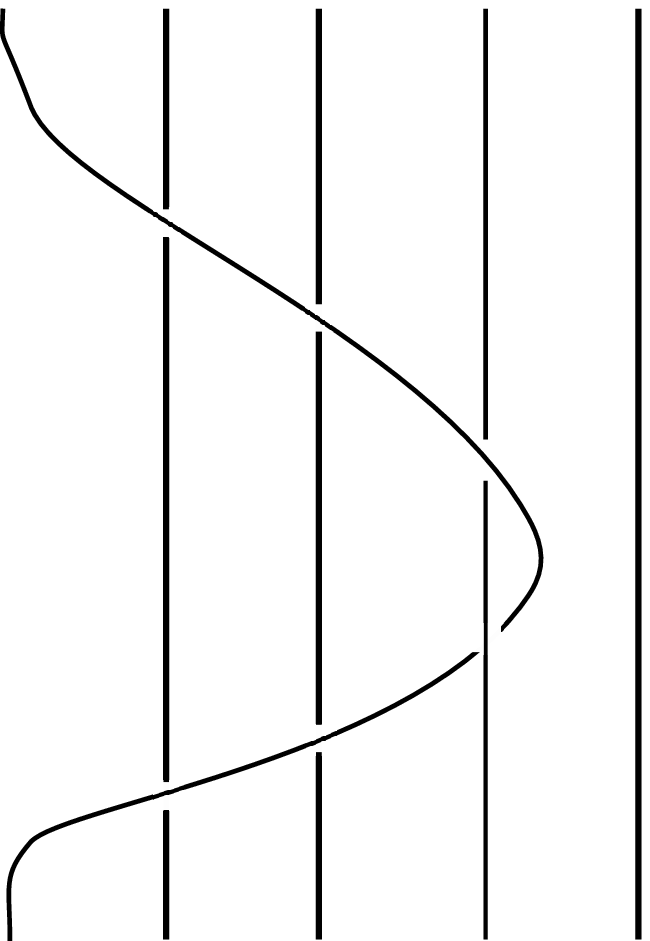}
\caption{$N=4$: Left, the path of $\alpha_3$;
 Right, the braid of $\alpha_3$}\label{alphafig}.
\end{center}
\end{figure}

To compute the entropy efficiency  we   need a 
specified set of generators   of $\cP_N$. For  $i = 1, \dots, N$, let 
\begin{equation}\label{alphadef}
\alpha_i = \sigma_1 \sigma_2 \dots \sigma_{i-1} \sigma_i^2 
\sigma_{i-1}\I \dots \sigma_2\I \sigma_1\I,
\end{equation}
where the $\sigma_i$ are the standard generators of $B_{N+1}$.
Note that we have adapted the convention that the stirrer $P$ is the 
left-most puncture and the obstacles are the next $N$
punctures.  Thus in terms of a stirring motion, the generator
$\alpha_i$ corresponds to the stirrer 
passing in front of the first $(i-1)$-obstacles from the left,
going around the $i^{th}$ obstacle (which is the $i+1^{st}$ 
puncture) once counterclockwise, and then returning to its initial position
passing once again in front of the first  $(i-1)$-obstacles
(see Figure~\ref{alphafig}).
The generators $\alpha_j$ are the members of the standard
generating set of the pure braid group $P_{N+1}$ in which
only the left-most strand is nontrivial. 
The next proposition follows from the theorem of Kra \cite{kra} and
Birman \cite{birmanbook} described above.
\begin{proposition}\label{pionefree}
The group $\cP_N$ of  \posp s with $N$ obstacles
is a free
group on the generators $\{ \alpha_1, \dots, \alpha_N\}$ and
is naturally isomorphic to $\pi_1(D_N, P)$.
\end{proposition}

\subsection{An adapted set of generators of $\pi_1(D_{N+1})$}\label{gensec}
\begin{figure}[ht]
\begin{center}
\includegraphics[width=.35\linewidth,angle=0]{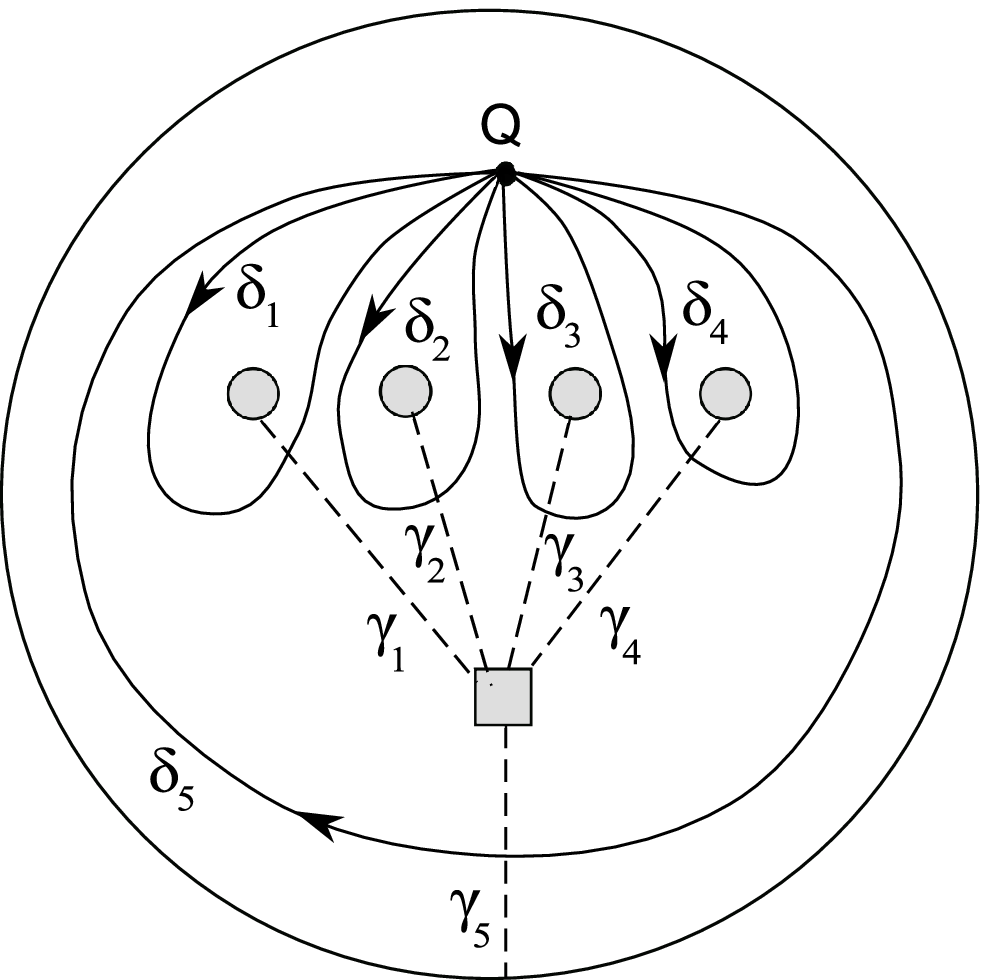}
\hskip 0.5in
\includegraphics[width=.35\linewidth,angle=0]{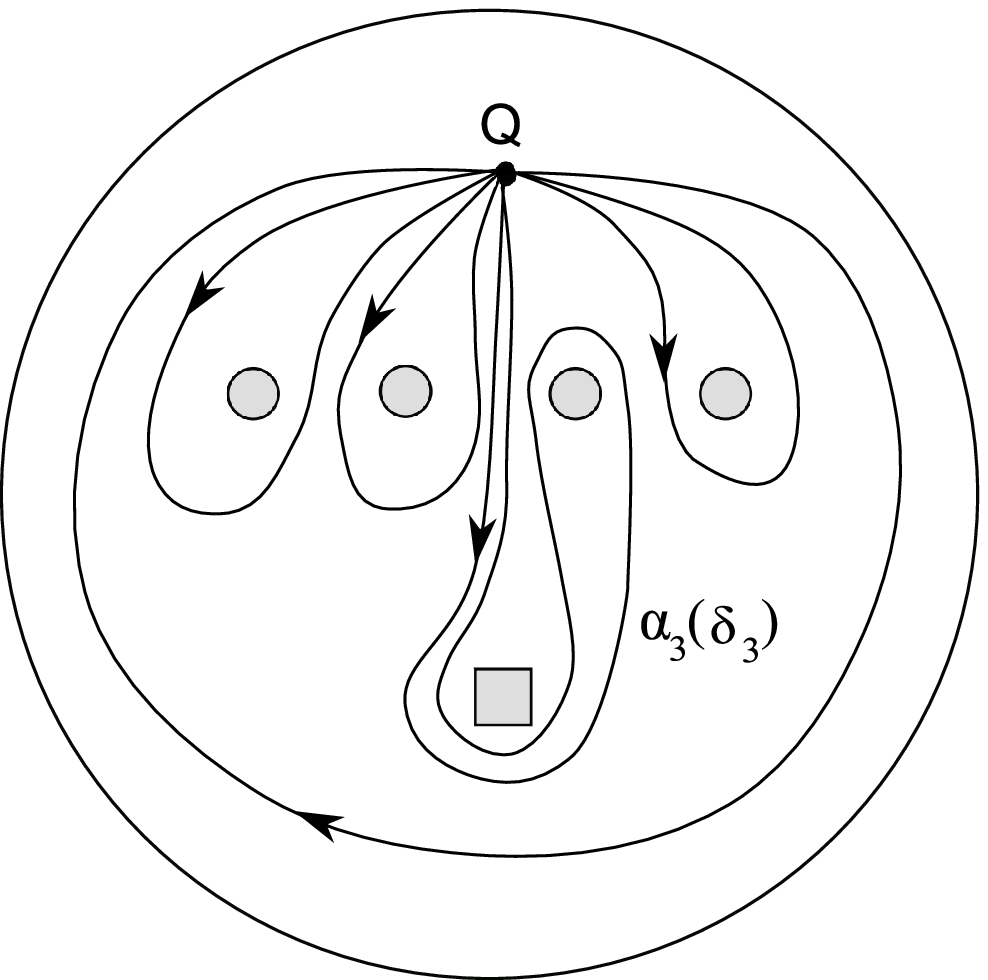}
\caption{Left, the adapted generators $\delta_j$ 
for $\pi_1 (D_5)$. The arcs $\gamma_j$ are used in the construction
of the covering spaces in \S\ref{coversec}. Right, the
action of $\alpha_3$.}
\label{genfig}
\end{center}
\end{figure}

In working with \posp s with $N$
obstacles it will be frequently
useful to use a collection of loops in $\pi_1(D_{N+1})$
which generate and are adapted to the protocols in  the sense that the
generators $\alpha_j$ of $\cP_N$ act particularly simply on them. 
The case
$N=4$ is shown in Figure~\ref{genfig}. The square represents the
stirrer $P$. Fix a basepoint $Q$, and for $i = 1, \dots, N$,
the loop $\delta_i$   begins at $Q$ and goes 
counterclockwise around the $i^{th}$ obstacle. The last
generator $\delta_{N+1}$ encloses the stirrer
and  all the obstacles going \textit{clockwise} around the disk.

The reason for the last convention is to get cleaner
 formulas for the action of the
$\alpha_j$ on the $\delta_i$.
We may assume that the generators $\{\alpha_1, \dots, \alpha_N\}$ 
of $\cP_N$ act while fixing the basepoint $Q$. Figure~\ref{genfig} right
shows the action of $\alpha_3$ when $N=4$.  Writing 
$\alpha_i(\delta_j) $ for the action of $\alpha_i$ 
on the loop $\delta_i$, we have
\begin{equation}\label{alphaact}
\begin{split}
\alpha_j(\delta_j) &= 
(\delta_{j-1}\I \delta_{j-2}\I \dots \delta_{1}\I\delta_{N+1}\I 
\delta_N\I  \dots \delta_{j+1}\I) 
\, \delta_j\,  
(\delta_{j+1} \dots  \delta_N \delta_{N+1}
\delta_{1} \dots \delta_{j-2} \delta_{j-1})  \\
\alpha_i(\delta_j) &= \delta_j \ \ \ \ \ \ \ \  \text{for} \ j \not= i. 
\end{split} 
\end{equation}

\subsection{Topological efficiency of \posp s} 
For a \posp\ $\chi\in\cP_N$, pick a \homeo\ $h\in\chi$ with
$h(Q) = Q$, and let $\chi_\sharp$ be the induced 
automorphism  of $\pi_1(D_{N+1}, Q) \isomorphic F_{N+1}$. Now 
$\eff(\chi_\sharp;\{\alpha_1\pmo, \dots, \alpha_N\pmo\} )$ is independent
of the choice of $h$ and $Q$, and    we  use the notation
\begin{equation*}
\eff(\chi) := \eff(\chi_\sharp;\{\alpha_1\pmo, \dots, \alpha_N\pmo\} ).
\end{equation*}
The maximal efficiency over the group $\cP_N$ of 
\posp s with $N$ obstacles is denoted here as 
\begin{equation*}
\Eff(N) := \sup\{\eff(\chi) : \chi\in\cP_N\}.
\end{equation*}

\begin{remark} To avoid potential confusion we emphasize that
a \posp\ with $N$ obstacles is determined by an element
of   the fundamental group of the $N$-punctured disk, $\pi_1(D_N)$,
but after also puncturing at the stirrer $P$ it yields
a mapping class on the $(N+1)$-punctured disk or a braid
on $(N+1)$-strings.
\end{remark}

\section{The lower bound}\label{lowerboundsec}
Since we are seeking $\Eff(N)$, the maximal efficiency of
 a \posp\ with $N$ obstacles, the
efficiency of any single $\chi\in\cP_N$ gives a lower bound. In the
first 
subsection we define  $\HSP_N\in\cP_N$ which numerically appears
to give the maximal efficiency, and we shall see  that
its topological efficiency converges as $N\raw \infty$
 to the upper bound given in \S\ref{upperboundsec}. Most  this section 
is concerned  with the computation of the topological entropy of $\HSP_N$.

\subsection{The hypertrophied stirring protocol}\label{hspdefinesec}
\begin{figure}[ht]
\begin{center}
\includegraphics[width=.4\linewidth,angle=0]{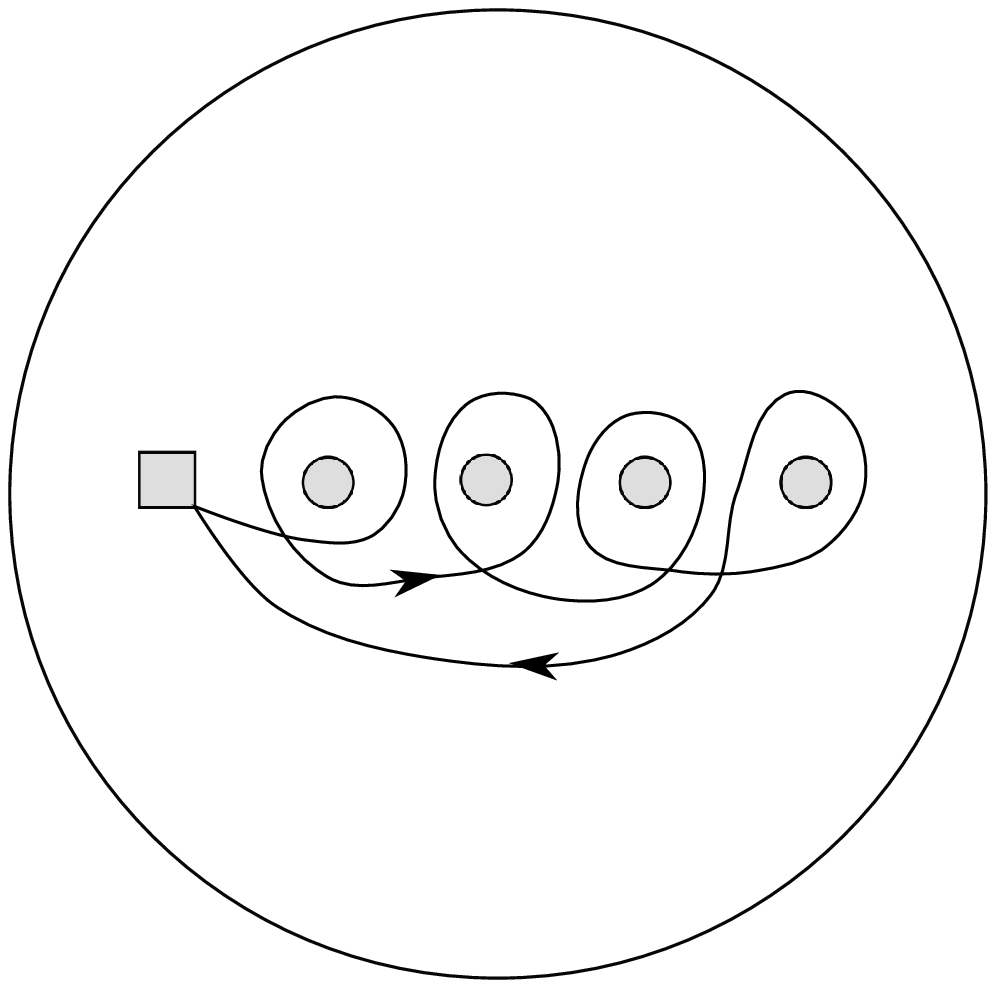}
\hskip 0.5in
\includegraphics[width=.4\linewidth,angle=0]{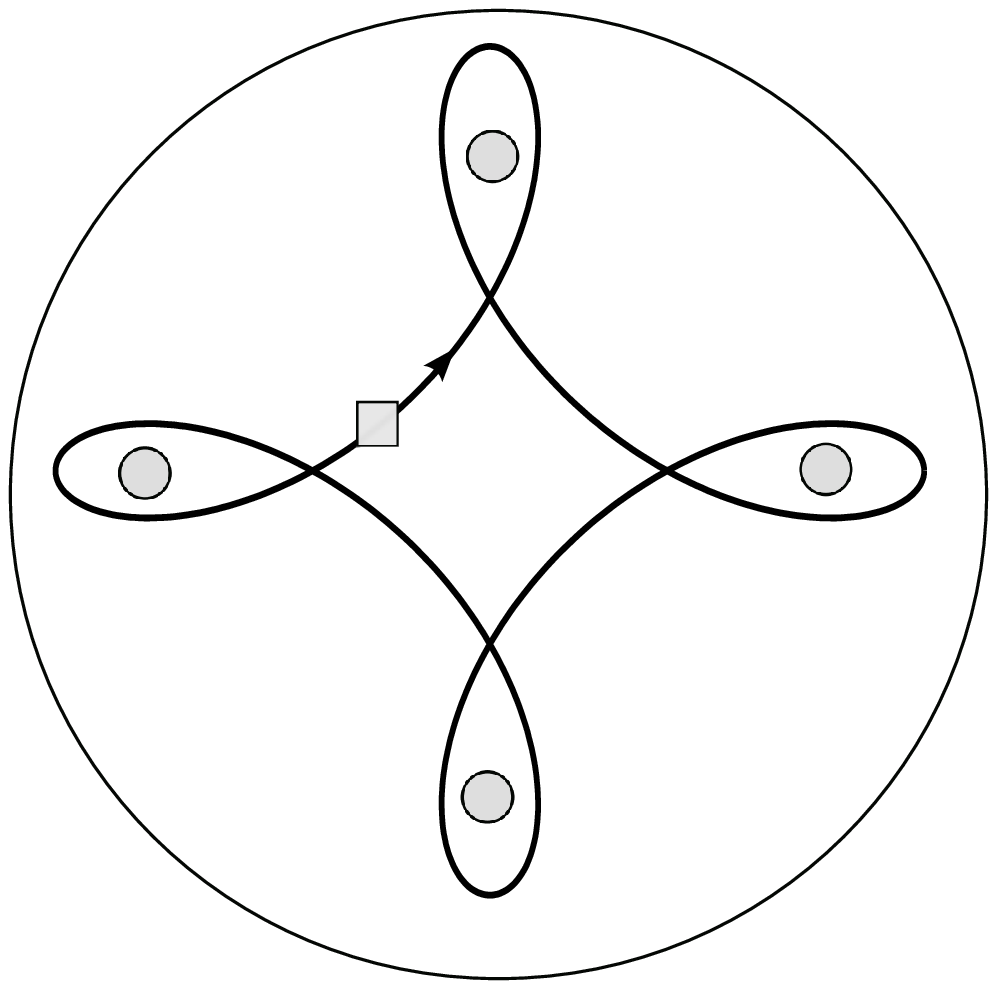}
\caption{ $\HSP_4$: Left, in
the braid configuration; Right, on a  hypotrochoid}
\label{hspfig}
\end{center}
\end{figure}

Using the generators $\alpha_i$ of $\cP_N$ we define 
the \de{hypotrochoid stirring protocol} with $N$ obstacles
as 
\begin{equation*}
\HSP_N := \alpha_1 \alpha_2 \dots \alpha_N.
\end{equation*}
As with any braid (see \S\ref{braidsec}),
in various contexts we treat  $\HSP_N$
as a physical braid,  a motion of a collection of points
(in this case  just the single stirrer is moving), and
a mapping class. Figure~\ref{hspfig} left shows the 
corresponding motion of the stirrer when $N=4$. Figure~\ref{hspfig} right
gives a topologically equivalent but more symmetric view
in which the path of the stirrer is a  hypotrochoid. 

\begin{figure}[ht]
\begin{center}
\includegraphics[width=.30\linewidth,angle=0]{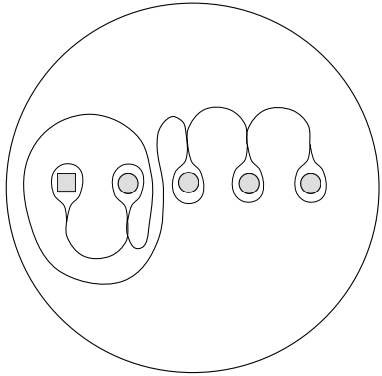}
\caption{The invariant train track for $\HSP_4$.}
\label{trackfig}
\end{center}
\end{figure}

The train track of a mapping class is a collapsed version the
unstable foliation (see \cite{penner}, \cite{casson}, \cite{primer}). Given
combinatorial data about a mapping class, it is always
computable in finite time by the Bestvina-Handel 
algorithm (\cite{BH}). The edges of the train track give a Markov partition
for a \pA\ representative $\phi$ in the class, and the induced
transition on the edges from the action of $\phi$ gives the
corresponding transition matrix.
Figure~\ref{trackfig} shows the invariant train track for 
 $\HSP_4$. Its invariance under the action of $\HSP_4$
as well as the fact that its transition matrix is strictly
positive can be confirmed by a rather tedious and long graphical
calculation or with a computer implementation
of the Bestvina-Handel algorithm (\cite{hall}). For other values 
of $N$ the repeated regularity of $\HSP_N$ ensures that
the train tracks in those cases is constructed by a repetition
of the simple one-prongs connected by a single edge seen
on the last three right-hand  punctures of $\HSP_4$.
The track demonstrates the important fact   
that the only singularities of the invariant foliations of
$\HSP_N$ are one-prongs at the punctures and 
a singularity on the outer boundary. Also 
note that the train track and its transition matrix prove that
each $\HSP_N$ is a \pA\ mapping class, though this is verifiable
by elementary means or from Kra's Theorem (\cite{kra}).

\subsection{The method for computing $h(\HSP_N)$}
In principle for any $N$ one can use the
invariant  train track  of $\HSP_N$  to  compute 
a Markov transition matrix for $\phi$. The dilation of
$\HSP_N$ is then the spectral radius of this matrix. 
However, this computation appears intractable  
and so we resort to homological methods.

As noted in \S\ref{thurstonsec}, the topological entropy of a \pA\
map $\phi$ can be computed as  the exponential
growth rate of words in the fundamental group
under the induced action of $\phi$. In practise,
this computation is very difficult. The corresponding
computation in homology is much simpler, but usually
just gives a lower bound (\cite{manning}).  Fried in \cite{friedtwisted}
develops the idea that the action of the induced
map on homology when lifted to finite covers
(or equivalently, on twisted homology) is usually computable,
often gives a better lower bound,  and in certain
cases, the bound is sharp.

 To implement this program,  
rather than working with individual finite covers 
it is simpler to construct a $\Z$-covering space $\tM$
and study its finite quotients.
The  first homology of $\tM$ is treated 
as a module over $\Z[t^{\pm 1}]$ with 
$t$ representing the generator of the deck group. 
If the given map $f$ (or the
 mapping class) lifts to $\tM$
and commutes with the deck transformation $T$, 
 then the lifted action of    $f$ 
is given by a matrix  $B(t) \in \Mat(m,  \Z[t, t\I])$,
where $m$ is the rank of $H_1(\tM)$ as a
$\Z[t^{\pm 1}]$ module. One may then show
(see Lemma~3.2  in \cite{bandboyland}) that  
a substitution of $t = \exp(2 \pi i p/q)$ into $B(t)$
gives the action
of $f$ on a piece of the first homology of
the $\Z_q$-cover constructed by taking the quotient of $\tM$
by $T^q$. 

In the case of the $n$-punctured disk, there is
a natural $\Z$-cover to which every mapping
class lifts and commutes with the deck transformations.
The matrix $B(t)$ giving the
action of the mapping class (or braid)  is the Burau 
representation and indeed, this is now the most
common descriptions of the topological meaning
of the  Burau representation.
Computing entropy of braids using substitutions
into the Burau representation is considered in
\cite{bandboyland}, and that paper is a good source for
more information on the techniques used in this
section.


\subsection{The $\Z$-covering spaces}\label{coversec}
The mapping classes representing \posp s are
rather special and so the  general cover used in 
Burau representation is not the best way to proceed.
Instead we construct
a pair of $\Z$-covering spaces,  $\tWN$ and $\tWNp$,
which are adapted to the particular nature of the action of \posp.
Both  $\tWN$ and $\tWNp$  
 are covering spaces of the $(N+1)$-punctured disk,
$D_{N+1}$. The constructions are   based on the 
generating loops $\{\delta_j\}$  for $\pi_1(D_{N+1})$ given in 
\S\ref{gensec}. We shall work with  simplicial one-chains 
 and give the relation to homology in 
Proposition~\ref{doublecoverok}.

\begin{figure}[ht]
\begin{center}
\includegraphics[width=\linewidth,angle=0]{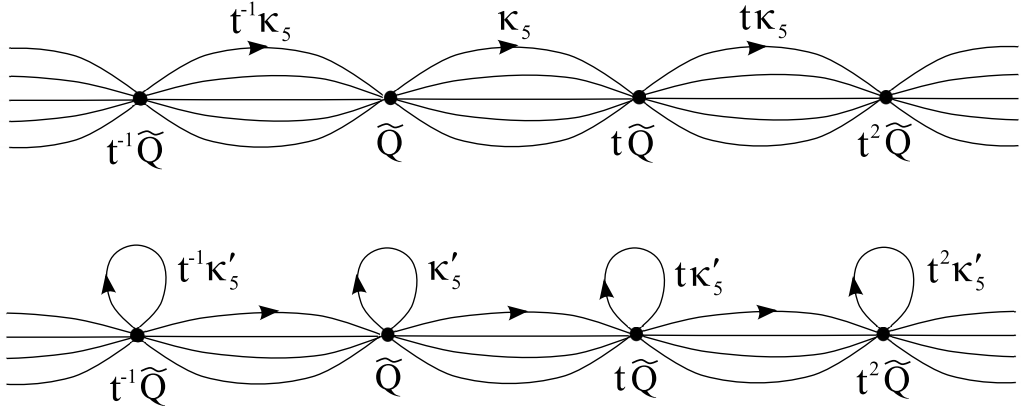}
\caption{Top, the spine of the covering
space $\tW_4$; bottom, the spine of the covering
space $\tWp_4$. Only selected simplices are labeled
and have their orientations shown.}
\label{coverfig}
\end{center}
\end{figure}

Geometrically,  both $\Z$-covering spaces are constructed by 
cutting along arcs ``dual'' to the generating loops
$\delta_1, \dots, \delta_{N+1}$ of $\pi_1(D_{N+1}, Q)$.
The arcs $\gamma_j$ for $j = 1, \dots, N$ connect the
$j^{th}$ obstacle to the stirrer, and $\gamma_{N+1}$ connects
the stirrer to the outer boundary (see Figure~\ref{genfig}).
 For the cover $\tWN$ one cuts along all these arcs and 
for $\tWNp$ one cuts along all of the arcs except $\gamma_{N+1}$.
After cutting along the arcs, copies of resulting  space are then
stacked and pasted into a $\Z$-cover.

For algebraic descriptions, recall that a covering space
with Abelian deck group of a manifold $M$ 
corresponds to a subgroup $H\subset H_1(M;\Z)$,
and the deck group of the resulting cover is  $H_1(M;\Z)/H$.
Thus a $\Z$-covering space of $D_{N+1}$ is given by an epimorphism
$\iota: H_1(D_{N+1};\Z) \raw \Z$. To define this epimorphism
for our two covering spaces, for $j = 1, \dots, N+1$,
we let $\Delta_j$  be the
element of  $H_1(D_{N+1};\Z)$ corresponding to the generating
loop $\delta_j$. Thus the set 
$\{\Delta_1, \dots, \Delta_{N+1}\}$ gives a basis
for $H_1(D_{N+1};\Z)$.
The covering spaces $\tWN$ and $\tWNp$ are then respectively 
defined by epimorphisms
\begin{equation}\label{iotadef}
\iota( \sum_{i=1}^{N+1} c_i \Delta_i) =  \sum_{i=1}^{N+1} c_i \
\ \ \text{and} \ \ \
 \iota'( \sum_{i=1}^{N+1} c_i \Delta_i) =  \sum_{i=1}^{N} c_i.
\end{equation}
In both cases the 
generator of the deck group  will be denoted $T$.

We now give the description of the one-chains in
$\tWN$ as a free $\Z[t^{\pm1}]$ module.
The wedge of $(N+1)$-circles given by the union of the generating
loops $\{\delta_j\}$ of $\pi_1(D_{N+1}, Q)$ is
denoted $\spine_N$ and it is a spine of $D_{N+1}$.
We work with the simplicial homology of
the one-skeleton $\spine_N$ and its lift to $\tWN$, denoted
$\tspine_N$. The  simplicial one-chains over $\Z$ 
are denoted $C_1(\spine_N)$ and $C_1(\tspine_N)$,
respectively.  Fix a lift  $\tQ$ of
the basepoint $Q$. For $j=1, \dots, N+1$, let
$\kappa_j$ be the oriented arc (or one-chain) in $\tWN$ which
is a lift of $\delta_i$ and connects $\tQ$ to $T(\tQ)$.
The collection $\{\kappa_1, \dots, \kappa_{N+1}\}$  
will be our basis for $C_1(\tspine_N)$, and
for $n\in\Z$, the formal expression $t^n \kappa_j$ refers
 to the one-chain $T_*^n(\kappa_i)$,
where $T_*: C_1(\tspine_N; \Z)\raw C_1(\tspine_N; \Z)$ is the
induced action of the deck transformation $T$. 
The one-chains in  $\tspine_N$ are thus the module $C_1(\tWN; \Z)
\isomorphic \Z[t^{\pm1}]^{N+1}$,
where the vector of Laurent polynomials
$(q_1, \dots, q_{N+1})$ with each
$q_j(t)\in \Z[t, t\I]$ is identified with the one-chain in 
$\tspine_N$ given by $\sum q_j(T_*) \kappa_j$.

The  one-chains in the spine $\tspine_N'$ of $\tWNp$
have a similar description. However in this case
$\delta_{N+1}$ doesn't lift to an arc but
rather lifts to a loop based at $\tQ$ which
is denoted $\kappa_{N+1}'$. We then use the simplicial one-chains
$\{\kappa_1, \dots, \kappa_{N}, \kappa_{N+1}'\} $
as a basis for $C_1(\tspine_N')$ over $\Z[t^{\pm1}]$ and proceed 
as with $\tWN$.

\subsection{The lifted action of a \posp}
We next consider the lifted action of a \posp\  $\chi\in\cP_N$.
Pick a \homeo\ $h\in\chi$ which fixes the basepoint $Q$
and let $\hDelta$ be its induced action on the
spine $\spine_N$. More precisely, let $\hDelta:\spine_N\raw \spine_N$
be a simplicial map homotopic to $r\circ h_{\vert \spine_N}$,
where $r:D_{N+1} \raw \spine_N$ is a  deformation
retract.
Since $h$ fixes all the punctures in  $D_{N+1}$, 
it acts as the identity on $H_1(D_{N+1};\Z)$. Thus $h$ lifts to    
$\tWN$ and commutes with all deck transformations. Let $\th$ be
the lift with $\th(\tQ) = \tQ$. The simplicial map
$\hDelta$ also lifts to $\thDelta:\tspine_N\raw \tspine_N$.  
Since $\thDelta$ commutes with the deck transformation, 
$\thDelta_*:C_1(\tspine_N)\raw C_1(\tspine_N)$ can be represented by 
a matrix $C(t)\in \Mat(N+1,  \Z[t^{\pm 1}])$.  The assignment 
$\chi\mapsto C(t)$ is a    
representation of the group of  \posp s with $N$ obstacles which we 
denote  $\Phi:\cP_N \raw \Mat(N+1,  \Z[t, t\I])$. 
It is clear that the construction does not depend on the
choice of $h\in\chi$ or its retract $\hDelta$ as long as
$\thDelta(\tQ) = \tQ$. If $\thDelta(\tQ) = T^n(\tQ)$,
then the corresponding matrix is multiplied by $t^n$.
Thus we adopt the notation $\tchi_* = \thDelta_*$ with
the understanding that it is always computed using
an $\thDelta$ with $\thDelta(\tQ) = \tQ$.
The same considerations for the cover $\tWNp$ yields
a representation  $\Phi':\cP_N \raw \Mat(N+1,  \Z[t^{\pm 1}])$,
and we sometimes write  $\Phi'(\chi) = C'(t)$ as an
abbreviated notation. 

Using their action on the generators $\delta_j$ it is straightforward
to now compute the action of the lift of generators $\alpha_j$ 
of $\cP_{N+1}$ on the chosen  basis
of simplicial one-chains in the $\Z$-covering space. 
 This procedure is described in \cite{friedtwisted}.
We let  $\talpha_j(\kappa_j)$ denote the lifted action of 
$\alpha_j$ on the simplex $\kappa_j$. 

\begin{example}
 We illustrate the general result by computing the lifted
action of $\alpha_3$ on $\tspine_4$. The action on $\pi_1$
from \eqref{alphaact} is 
\begin{equation}\label{alphathree}
\alpha_3(\delta_3) = \delta_2\I \delta_1\I \delta_5\I \delta_4\I
\delta_3 \delta_4 \delta_5 \delta_1 \delta_2,
\end{equation}
and $\alpha_3(\delta_j) = \delta_j$ for $j\not = 3$. For the
lifted action to the generators of the simplicial chains 
in $\tspine$, it is immediate that $\talpha_3(\kappa_j) =
\kappa_j$ for $j\not = 3$, so we need only compute 
 $\talpha_3(\kappa_3)$. Since our convention is to
represent mapping classes with \homeos\ whose
lifts fix $\tQ$, the image of $\kappa_3$ must be an
arc beginning at $\tQ$. Since from \eqref{alphathree}
$ \alpha_3(\delta_3)$ begins with $\delta_2\I$,
the beginning of $\talpha_3(\kappa_3)$ must
traverse a copy of $\kappa_2$ in the negative direction
and thus must in fact be $-t\I \kappa_2$. The image of
$\delta_3$ continues with $\delta_1\I$ and so the
image of $\kappa_3$ continues with $-t^{-2} \kappa_1$.
Continuing in this manner the next stage of the image
 $\talpha_3(\kappa_3)$ is always forced by the next 
generator (or its inverse) in the image $ \alpha_3(\delta_3)$.
The full image is then
\begin{equation*}
\begin{split}
\talpha_3(\kappa_3) &= -t\I \kappa_2 -t^{-2} \kappa_1 -t^{-3} \kappa_5
-t^{-4} \kappa_4 + t^{-4} \kappa_3 + 
t^{-3} \kappa_4 + t^{-2} \kappa_5
+ t^{-1} \kappa_1 + \kappa_2\\
&= \kappa_2 + t\I(\kappa_1 - \kappa_2) + t^{-2}(\kappa_5 - \kappa_1)
+ t^{-3}(\kappa_4 - \kappa_5) + t^{-4}(\kappa_3 -\kappa_4)\\
\end{split}
\end{equation*}
The process of deriving the
formulas for the action of $\alpha_3$ lifted to
$\tspine'_4$ is similar, but now instead of traversing
the arcs corresponding to lifts of $\delta_5$ to $\tspine_4$, one
  goes around the loops $t^n \kappa_5'$ which are the lifts
of $\delta_5$ to $\tspine'_4$.
\end{example}

The same procedure in general yields the following.
\begin{equation}\label{talphaact}
\begin{split}
\talpha_1(\kappa_1) &= 
\kappa_{N+1} + t^{-1}(\kappa_{N} - \kappa_{N+1}) + \dots + 
 t^{-N}(\kappa_{1} - \kappa_{2})\\
\talpha_j(\kappa_j) &= 
\kappa_{j-1} + t^{-1}(\kappa_{j-2} - \kappa_{j-1}) + \dots + 
 t^{-(j-2)}(\kappa_{1} - \kappa_{2}) + 
 t^{-(j-1)}(\kappa_{N+1} - \kappa_{1}) + \\ 
 &\ \ \ \ t^{-j}(\kappa_{N} - \kappa_{N+1}) + 
 t^{-(j+1)}(\kappa_{N-1} - \kappa_{N}) + \dots 
 t^{-N}(\kappa_{j} - \kappa_{j+1}) \ \  \text{when}\ j > 1 \\ 
\talpha_i(\kappa_j) &= \kappa_j \hskip 60pt \text{when}\ j\not= i
\end{split}
\end{equation}
Letting   $\talpha'_i$ the lifted action of 
$\alpha_i$ on one-chains in $\tspine'$ we have in general 
\begin{equation}\label{talphapact}
\begin{split}
\talpha'_1(\kappa_1) &= 
\kappa_1 +  \kappa_{N+1} +  t^{-1}(\kappa_{N} - \kappa_1- \kappa_{N+1}) + 
t^{-2}(\kappa_{N-1} - \kappa_{N}) + \dots + 
 t^{-(N-1)}(\kappa_{2} - \kappa_{3})\\
\talpha'_2(\kappa_2) &= 
t \kappa_{N+1} +  (\kappa_{N} - \kappa_{N+1}) + 
t^{-1}(\kappa_{N-1} - \kappa_{N}) + \dots + 
 t^{-(N-1)}(\kappa_{1} - \kappa_{2})\\
\talpha'_j(\kappa_j) &= 
\kappa_{j-1} + t^{-1}(\kappa_{j-2} - \kappa_{j-1}) + \dots + 
 t^{-(j-3)}(\kappa_{2} - \kappa_{3}) + 
 t^{-(j-2)}(\kappa_{1} - \kappa_{2} + \kappa_{N+1}) + \\
&\ \ \ \ \  t^{-(j-1)}(\kappa_{N} - \kappa_{1}- \kappa_{N+1}) + 
  t^{-j}(\kappa_{N-1} - \kappa_{N}) + 
    \dots 
+ t^{-(N-1)}(\kappa_{j} - \kappa_{j+1})  \  \text{when}\ j > 1\\  
\talpha'_i(\kappa_j) &= \kappa_j \hskip 60pt \text{when}\ j\not= i
\end{split}
\end{equation}

\subsection{The double cover}\label{doublecoversec}
The two-fold covering spaces
formed by taking the quotient of $\tWN$ and  $\tWNp$ by
the square of the deck transformation will be central
to our calculations. We recall some basic facts
concerning the connection of  various objects
in the full $\Z$-covering spaces, the base space
and this two-fold cover, and then compute the action
on homology in this cover. We will initially just discuss
$\tWN$ and comment on  $\tWNp$ afterwards.
To simplify notation fix $N$ and drop it from the notation,
and so now $\tW = \tWN$.

 Recall that the $\Z$-covering space
$p:\tW\raw D_{N+1}$  corresponds to an epimorphism $\iota:H_1(D_{N+1};\Z) 
\raw \Z$. The two-fold covering space 
$\pt:\tWt\raw D_{N+1}$ corresponds to the epimorphism
$\iotat:H_1(D_{N+1};\Z) \raw \Z_2$ where $\iotat := \pi^{(2)} \circ
\iota$, with  $\pi^{(2)}:\Z\raw \Z_2 = \Z/(2\Z)$ the projection.
Thus $\tWt = \tW/T^2$, and we denote the deck transformation
of  $\tWt$ as $\Tt$. The corresponding one-skeleton of $\tWt$  is
denoted $\tspinet$.

If $\tau$ is an indeterminate 
with $\tau^2 = 1$,  then $C_1(\tspinet)$ is a $\Z[\tau]$
module with  $C_1(\tspinet)\isomorphic \Z[\tau]^{N+1}$. 
The action of $\talpha_i$ is denoted as $\talphat_i$ and
using \eqref{talphaact} we get 
\begin{equation}\label{tauaction}
\begin{split}
\talphat_j(\kappa_j) &=    (-1)^{j} (1-\tau) \kappa_1 + 
  (-1)^{j+1} (1-\tau) \kappa_2 + \dots + 
 (1-\tau) \kappa_{j-1} +  
   \kappa_{j} - \\ 
& \ \ \ \ \ (1-\tau) \kappa_{j+1} + \dots +
 (-1)^{(j+1)} (1-\tau)  \kappa_{N+1}\ \ \text{when}\ N \ \text{is even}\\
\talphat_j(\kappa_j) &=     (-1)^{j} (1-\tau) \kappa_1 + 
  (-1)^{j+1} (1-\tau) \kappa_2 + \dots + 
 (1-\tau) \kappa_{j-1} - 
  \tau \kappa_{j} + \\ 
& \ \ \ \ \  (1-\tau) \kappa_{j+1} + \dots +
 (-1)^{(j+1)} (1-\tau)  \kappa_{N+1}\ \ \text{when}\ N \ \text{is odd}\\
\talphat_i(\kappa_j) &= \kappa_j \ \ \ \ \
 \text{for all}\ N\ \text{when}\ j\not= i
\end{split}
\end{equation}
 Similarly, $\tWpt = \tWp/T^2$ 
is defined by the epimorphism $\iotat:= \pi^{(2)} \circ
\iota'$. Its one-skeleton is $\tspinept$
 and  $C_1(\tspinept)\isomorphic \Z[\tau]^{N+1}$
and using \eqref{talphapact} the action of $\talpha_i$ is 
\begin{equation}\label{taupaction}
\begin{split}
\talphapt_j(\kappa_j) &=     (-1)^{j} (1-\tau) \kappa_1 + 
  (-1)^{j+1} (1-\tau) \kappa_2 + \dots + 
 (1-\tau) \kappa_{j-1} - 
  \tau \kappa_{j} +  \\
& \ \ \ \ \ (1-\tau) \kappa_{j+1} + \dots +
 (-1)^{j} (1-\tau)  \kappa_{N+1}\ \ \text{when}\ N \ \text{is even}\\
\talphapt_j(\kappa_j) &=     (-1)^{j} (1-\tau) \kappa_1 + 
  (-1)^{j+1} (1-\tau) \kappa_2 + \dots + 
 (1-\tau) \kappa_{j-1} + 
   \kappa_{j} -  \\
& \ \ \ \ \ (1-\tau) \kappa_{j+1} + \dots +
 (-1)^{j} (1-\tau)  \kappa_{N+1}\ \ \text{when}\ N \ \text{is odd}\\
\talphapt_i(\kappa_j) &= \kappa_j \ \ \ \ \
 \text{for all}\ N\ \text{when}\ j\not= i
\end{split}
\end{equation}

Let $\chi$ be a \posp\ treated as a mapping class. The projection
of the action $\tchi_*$ on $C_1(\tspine)$ to one 
on $C_1(\tspinet)$ is denoted $\tchit_*$. The action
of $\chi$ on $H_1(\tspinet;\Z) \isomorphic H_1(\tWt;\Z)$
is written  $\tchit_\star$. The same notations are used for 
$\tWpt$, but with a prime attached.

The next proposition gives the connection of the substitution
$t=-1$ into the matrices $C(t)$ and $C'(t)$ and
the growth rate of the action of $\chi$ on the first
homology of the two-fold cover $\tWt$. 
It  follows from Lemma $3.2$ in \cite{bandboyland} and
the fact that the growth on chains is the same as that
on homology.  For sake of having the paper 
self-contained, we  give the elementary proof.
\begin{proposition}\label{doublecoverok}
Assume that $\chi\in\cP_N$ is a given \posp. As constructed
above,  $\tWt$ and $\tWpt$ are double covers
and $\tchit_\star$ and $\tchipt_\star$ are the induced
actions on $H_1(\tWt)$ and $H_1(\tWpt)$, respectively,
and  $\Phi$ and $\Phi'$ are representations of
$\cP_N$.
Letting $C(t) = \Phi(\chi)$ and
$C'(t) = \Phi'(\chi)$,
 then
\begin{equation*}
\rho( \tchit_\star )= \rho( C(-1)) \ \
\text{and}\ \   \rho( \tchipt_\star )= \rho( C'(-1)).
\end{equation*}
\end{proposition}

\proof\ 
 We will initially just discuss
$\tWN$ and comment on  $\tWNp$ afterwards.

Since $(\Tt)^2 = \id$, the eigenvalues of
its action on simplicial one-chains in  $\tspinet$
 are $\pm 1$, and we denote the 
corresponding eigen-space decomposition as $C_1(\tspinet) = V_+ \osum V_-$.
 The subspaces $V_\pm$ are also characterized as 
$V_+ = \{ c + \Tt_*(c) : c \in C_1(\tspinet)\}$ and 
$V_- = \{ c - \Tt_*(c) : c \in C_1(\tspinet)\}$. 
Since $\tchit_*$ commutes with $\Tt_*$, we have that
 $\tchit_*$ preserves the decomposition
$C_1(\tspinet) = V_+ \osum V_-$.

Treating  $C_1(\tspinet)\isomorphic \Z[\tau]^{N+1}$, 
the eigen-subspaces   have alternative descriptions as  
$V_+ = \{ c + \tau c : c \in C_1(\tspinet)\}$ and 
$V_- = \{ c -\tau c : c \in C_1(\tspinet)\}$.  
The homomorphism $\tchit_*$ acts on $C_1(\tspinet)$
 as the matrix $C(\tau)$, and we will write  
 $C(\tau) = B_0 + \tau B_1$ with 
$B_i\in\Mat(N+1, \Z)$. Thus acting on an element
of $V_+$ we have $(c + \tau c) C(\tau) = 
(c B_0 + c B_1) + \tau(c B_1 + c B_0) = (c + \tau c) (B_0 + B_1)
= (c + \tau c) C(1)$
Similarly acting on an element of $V_-$ we have  
$(c - \tau c) C(\tau) = (c - \tau c) C(-1)$.
So  $\tchit_*$ acts on $V_+$ by $C(1)$ and on 
$V_-$ by $C(-1)$. Substituting $\tau=1$ into
\eqref{tauaction}, we find $C(1) = I$ (this also
follows from the fact that 
$\pt_*:  C_1(\tspinet) \raw C_1( D_{N+1})$ 
is an isomorphism when restricted to $V_+$).
Since $\det(C(-1)) = 1$, $\rho(C(-1)) \geq 1$ and so 
 $\rho(\tchit_*) =   \rho(C(-1))$.

To find out what this means for homology
we go back to working over $\Z$.
For $1\leq i \leq N+1$, let $\wkappa_i = \Tt_*(\kappa_i)$,
and so $\{\kappa_i, \wkappa_i\}$ forms a basis
for $C_1(\tspinet;\Z)$. Now
  for  $1\leq i \leq N$, let
$\Gamma_i = \kappa_i - \kappa_{i+1}$ and
$\wGamma_i = \wkappa_i - \wkappa_{i+1}$.
Letting $\hGamma = \kappa_{N+1} + \wkappa_{N+1}$
it is easy to check that 
$\{\Gamma_1,\dots, \Gamma_N, \wGamma_1,\dots, \wGamma_N, \hGamma,
\kappa_{N+1}\}$ is also a basis for $C_1(\tspinet;\Z)$
and we define $V$ as the linear span
 $V = \langle \Gamma_1,\dots, \Gamma_N, 
\wGamma_1,\dots, \wGamma_N, \hGamma \rangle$
Thus $C_1(\tspinet;\Z) = V \osum \langle \kappa_{N+1}\rangle$.
Now $\hGamma$ and   all of the $\Gamma_i$ and 
$\wGamma_i$ are cycles because there are no two cells in $\tspinet$.
 In fact, $V$ is   a concrete realization of the simplicial
homology $H_1(\tspinet;\Z)$. Either by direct
calculation using \eqref{tauaction}, or from general principles,
$\tchit_*(V) = V$. Also from \eqref{tauaction}, $\talphat_j(\kappa_{N+1})
= \kappa_{N+1}$, and so we see that the spectral radius of
 $\tchit_*$ acting on all of $C_1(\tspinet;\Z)$ is the same as the spectral
radius of its restriction to $V$. But $\tchit_*$ restricted 
to $V$ is exactly $\tchit_\star$, the induced action of $\tht$ on
$H_1(\tWt;\Z)$ and so by the conclusion of the  previous
paragraph, $\rho(\tchit_\star) = C(-1)$.

For the other cover $\tWpt$, the proof on the chain level that  
$\rho(\tchipt_*) =   \rho(C'(-1))$ is exactly the
same as for $\tWt$. A slight alteration is required 
when passing to homology because  $\delta_{N+1}$ lifts
to a loop in $\tWpt$: for  $1\leq i \leq N-1$, let
$\Gamma_i = \kappa_i - \kappa_{i+1}$ and
$\wGamma_i = \wkappa_i - \wkappa_{i+1}$, and finally  
let $\hGamma = \kappa_{1} + \wkappa_{1}$.
It is easy to check that 
$\{\Gamma_1,\dots, \Gamma_{N-1}, \kappa_{N+1},
 \wGamma_1,\dots, \wGamma_N, \kappa_{N+1}, \hGamma,
\kappa_{1}\}$ is  a basis for $C_1(\tspinept;\Z)$,
and we let
 $V' = \langle \Gamma_1,\dots, \Gamma_{N-1}, \kappa_{N+1},
 \wGamma_1,\dots, \wGamma_N, \kappa_{N+1}, \hGamma\rangle$. Thus
 $C_1(\tspinept;\Z) = V' \osum \langle \kappa_{1}\rangle$, and 
 $V'$ is a concrete realization of  $H_1(\tspinept;\Z)$ with 
$\tchipt_*(V') = V'$. From \eqref{taupaction},  
$\talphapt_1(\kappa_{1}) =\tau \kappa_1 + \dots$ when $N$ is even,
and $\talphapt_1(\kappa_{1}) =  \kappa_1 + \dots$ when $N$ is
odd, and
 $\talphapt_j(\kappa_{1}) = \kappa_1$ for   $j\not = 1$, and so 
the spectral radius of $\tchipt_*$
acting on all of $C_1(\tspinept;\Z)$ is the same as the spectral
radius of its restriction to $V'$. But recall
that   $\tchipt_*$ restricted 
to $V'$ is exactly $\tchipt_\star$, the induced action of $\tchipt$ on
$H_1(\tspinept;\Z)$. Therefore  by the conclusion of the first
paragraph, $\rho(\tchipt_\star) = C'(-1)$.
\QED\

\subsection{Another representation}\label{anotherrepsec}
Proposition~\ref{doublecoverok} shows the significance of 
the substitution  $t=-1$  into the matrices $C(t)$ and $C'(t)$.
It will turn out that the matrix of significance for computing
$h(\HSP_N)$ will be $C(t)$ when $N$ is even and $C'(t)$ when
$N$ is odd.
Conveniently, 
the substitution of $t=-1$ 
into  \eqref{talphaact} (or $\tau = -1$ into \eqref{tauaction})  
when $N$ is even
yields an identical formula
to the result of substituting $t=-1$ into
\eqref{talphapact}  when $N$ is odd.

Now note that for all $j$, $\talphat_j(\kappa_{N+1}) = \kappa_{N+1}$,
and so for  any product of the $\talphat_j$,
$\kappa_N$ is  an eigenvector with eigenvalue $1$.
Since we are interested in exponential growth rates,
we may ignore $\kappa_{N+1}$.  For $1\leq k \leq N$,
define
the $(N\times N)$-matrix $T^{(k)}$  by 
\begin{align}\label{repeq}
T^{(k)}_{kj} &= 2\cdot (-1)^{k - j +1} &\text{for} \ 1 \leq j \leq k-1\notag\\
T^{(k)}_{kk} &= 0 &\\
T^{(k)}_{kj} &= 2\cdot (-1)^{j-k} &\text{for} \  k+1 \leq j \leq N\notag\\
T^{(k)}_{ij} &= 0 &\text{for} \ i\not = k,\notag
\end{align}
and let 
\begin{equation*}
E^{(k)} = I + T^{(k)}.
\end{equation*}
  Extend the assignment
$\halpha_k\mapsto E^{(k)}$ to a homomorphism 
 $\Psi:\cP_N\raw \SL(N, \Z)$. Thus $\Psi(\chi)$ 
is the matrix $C(-1)$ with the last row and column
deleted for $N$ even and
$\Psi(\chi)$ 
is the matrix $C'(-1))$ with the last row and column
deleted for $N$ odd.

\subsection{A trace computation}
Fix an $N\geq 2$. Recall that
$\HSP_N = \alpha_1 \dots \alpha_N$. 
We will compute  $\trace(\Psi(\HSP_N))$ inductively
using the  initial segments of the product 
 $\alpha_1 \dots \alpha_k$ for $k = 1, \dots, N$.
These products have the matrix representation   
$H^{(k)} := \Psi(\alpha_1 \dots \alpha_k) = 
E^{(1)} E^{(2)} \dots E^{(k)}$.
\begin{lemma}\label{tracelem}
With $\Psi$ and $\HSP_N$ as defined above, 
 $\trace(\Psi(\HSP_N)) = -3^N + 3N +1$.
\end{lemma}

\proof\ It follows easily from the definitions that
for any $1\leq a,b \leq N$,
\begin{equation}\label{basicit}
\Hko_{a,b} = \Hk_{a,b} + \Hk_{a, k+1}\cdot T^{(k+1)}_{k+1, b}.
\end{equation}
We now prove by induction on $k < N$ that for $m = 1, \dots, N-k$,
\begin{align}
\Hk_{i, k+m} &= 2\cdot 3^{k-i} (-1)^{i + k + m} & 1\leq i \leq k\notag\\
\Hk_{i, k+m} &= 0 & k+1 \leq i \leq N\ \text{and}\ i\not= k+m
 \label{induct}\\
\Hk_{k+m, k+m} &= 1 &\ \notag  
\end{align}
The case $k=1$ is trivial. So assume 
\eqref{induct} is true for $k$. Thus using 
\eqref{basicit}. for  $1\leq i \leq k$ and  $m = 2, \dots, N-k$,
\begin{equation*}
\begin{split}
\Hko_{i, k+m} &= \Hk_{i, k+m}  + \Hk_{i, k+1} T^{(k+1)}_{k+1, k+m}\\
&=  2\cdot 3^{k-i} (-1)^{i + k + m} + 
 2\cdot 3^{k-i} (-1)^{i + k + 1}\cdot  2 (-1)^{k+m-(k+1)} \\
&= 2\cdot 3^{k+1-i} (-1)^{i + k + 1 + m},
\end{split}
\end{equation*}
as required. The other values of $i$ are special cases
easily checked, completing the proof of \eqref{induct}.

Now we prove that for $1 \leq k \leq N$, 
 $\trace(\Hk) = -3^k + 2k + N +1$, again by induction
on $k$, with the case $k=1$ being trivial. Using
\eqref{basicit} and the formulas \eqref{induct},
\begin{equation*}
\begin{split}
\trace(\Hko) &= \trace(\Hk) + \sum_{a=1}^N \Hk_{a, k+1} T^{(k+1)}_{k+1,a}\\
&= \trace(\Hk) + \sum_{a=1}^k 2\cdot 3^{k-a}  (-1)^{a + k +1} 2 (-1)^{k-a}\\
& = -3^{k+1} + 2(k+1) + N +1. 
\end{split}
\end{equation*}
The statement in the lemma now follows by letting $k=N$.
\QED\ 

\subsection{The orientation double cover of a \pA\ map}
A \pA\ map $\phi$ is characterized by the existence of 
a pair of transverse invariant
measured foliations. If these foliations are oriented, 
it is now standard
that the dilation of the \pA\ map is the spectral radius
of the action of $\phi$ on first homology.
 For \pA\ maps on the
punctured disk and in many other cases, 
a \pA\ map's foliations are not oriented. However, one  
 can lift $\phi$ and its foliations to
a new \pA\ map $\ophi$ on 
a branched cover $\oM$ where there foliations are oriented,
and thus the dilation is detectable using homology.
The construction of this cover requires
a fair amount of \textit{a priori} information about
$\phi$, and the exact connection between $\phi$
and $\ophi$ can be difficult to ascertain, so in 
many situations the orientation double cover is not a useful
tool.  However, because the family $\HSP_N$ lives on a genus zero
surface and its foliations are of a particularly simple 
type, it will turn out that its orientation double cover
always corresponds to one of  the two-fold covers
 finite covers constructed
by taking quotients of the $\Z$-covering spaces
$\tWt$ and $\tWpt$ which were considered in \S\ref{doublecoversec}. 

A few definitions are required to give a characterization
applicable to the case at hand.
Assume $\phi$ is a \pA\ map on the $n$-punctured disk $D_n$ 
with unstable foliation $\cF^u$. Let $S$ be the collection
of singularities of $\cF^u$ which are of odd order and
contained in the interior of $D_n$. In addition,   $P$ denotes
the collection of punctures of $D_n$ and $S'\subset P$
is the subset consisting of those punctures at which the singularity
of $\cF^u$ is of odd order. The orientation double
cover   of $\phi$ is the two-fold  cover branched
over $S\cup S'$. More formally, first further puncture 
$D_n$ at $S$ to form $M' = D_n - S$. For each $p_i\in P \cup S$,
 let $\Gamma_i $
be a small loop going counterclockwise around $p_i$, and so
the homology classes $\{[\Gamma_i]\}$ are a basis for
$H_1(M';\Z)$. If $I$ is the set of indices $i$ so that 
  $p_i\in S' \cup S$,   define the epimorphism
 $H_1(M';\Z)\raw \Z_2$  
by
\begin{equation*}
\sum c_i [\Gamma_i] \mapsto \sum_{i\in I} c_i \mod 2.
\end{equation*}
The corresponding two-fold cover over $M'$ is denoted
$\hM'$. Finally, fill in the punctures in $ \hM'$ 
that project to points in $S$ to form  $\hM''$ which
is a branched covering space over $D_n$.  
The following is standard: for example,
it  follows from Lemma $2.2$ in \cite{bandboyland}.
\begin{proposition}\label{odouble}
Let $\phi$ be \pA\ map on the $n$-punctured disk $D_n$ and
 $\hM''$ the branched double cover of $D_n$ formed by
branching over the singularities of the foliations of
$\phi$ which have odd order as was just constructed.
Then  $\hM''$ is the orientation double cover  of $\phi$. 
\end{proposition}

\subsection{The entropy of $\HSP_N$ and the lower bound}
We now have all the ingredients needed to obtain our
lower bound.
\begin{theorem}\label{lowerboundthm}
The topological entropy
of the hypotrochoid stirring protocol with $N$ obstacles is given
by the log of the spectral radius of the matrix
 $H^{(N)}$ defined in the previous subsection,
\begin{equation}\label{matcomp}
h(\HSP_N) = \log(\rho(H^{(N)})).
\end{equation}
As a consequence,
\begin{equation*}
  \log(\frac{3^N -3N - 1}{N}) \leq h(\HSP_N),
\end{equation*} 
and thus the maximal efficiency of a \posp\ with
$N$ obstacles satisfies
\begin{equation*}
 \frac{\log(3^N - 3N - 1) - \log(N)}{N} \leq \Eff(N).
\end{equation*} 
\end{theorem}

\proof\ The main observation we need is that the double covers
$\tWt$ and $\tWpt$ from \S\ref{doublecoversec}
 are identical to the
orientation double covers of a \pA\ map $\phi\in\HSP_N$ when $N$ is even
and odd, respectively. 
Once again we discuss the case of $\tWt$ first.
Recall that $\tWt$ corresponds to 
the epimorphism $\rhot:H_1(D_{N+1};\Z) \raw \Z_2$ given
by $\rhot(\sum_{i=1}^{N+1} c_i \Delta_i)\equiv 
\sum_{i=1}^{N+1} c_i  \mod 2$ 
with $\Delta_i$ the basis given in \S\ref{coversec}.
 Let $\Omega$ be the homology 
class of a small loop going once counterclockwise around
the stirrer $P$. We then have $\Omega = -\sum^{N+1}_{i=1} \Delta_i$ and
so $\rhot(\Omega) \equiv 1 \mod 2$ when $N$ is even. As noted
in \S\ref{hspdefinesec}, the foliations
of the \pA\ representative  $\phi\in\HSP_N$ have their only
singularities at the punctures and outside boundary, and the
singularities at the punctures all have odd order. Thus  by
Proposition~\ref{odouble}, the orientation double cover of $\phi$ is
defined by the epimorphism,
$\orho(\Delta_i) \equiv 1\mod 2$ for $i=1, \dots, N$ and 
$\orho(\Omega) \equiv 1 \mod 2$. Thus the orientation double cover
of $\HSP_N$ is exactly $\tWt$ when $N$ is even.

Again for even $N$, letting $C(t) = \Phi(HSP_N)$, 
 Proposition~\ref{doublecoverok} says that  
$\rho(\tphit_\star) = \rho(C(-1))$, where $\tphit$ is
the lift of the \pA $\phi\in\HSP_N$ to $\tWt$.
But since $\tWt$ is also the orientation
double cover for $\phi$, we have $h(\phit) = 
\log(\rho( \tphit_\star))$. Since entropy is preserved by finite
covers, we also have $h(\phi) = \log(\rho( \tphit_\star))$. Finally,
as remarked in \S\ref{anotherrepsec},
 $\rho(H^{(N)}) = \rho( C(-1))$, finishing
the proof when $N$ is even.

To obtain the result when $N$ is odd, we use the
two-fold cover $\tWpt$. It  corresponds to the epimorphism
$\rhopt(\sum^{N+1}_{i=1} c_i \Delta_i) \equiv \sum^{N}_{i=1} c_i  \mod 2$.
Thus when $N$ is odd, $\rhopt (\Omega) \equiv 1 \mod 2$, and once
again $\tWpt$ is the orientation double cover of the \pA\ 
$\phi\in \HSP_N$. The rest of the argument yielding
\eqref{matcomp} is the same as when $N$ is even.

Since for any $(N\times N)$-matrix $C$, 
$\rho(C) \geq \trace(C)/N$, the last two statements of the theorem
follow from Lemma~\ref{tracelem}.
\QED\ 

\subsection{The case $N=2$}\label{casetwo}
The case of two obstacles requires separate 
consideration. While Theorem~\ref{lowerboundthm} is
true for $N=2$, the lower bound is simply 
$0 \leq h(\HSP_2)$ and indeed, it is easy to
check that $\rho(\HSP_2) = 1$. Fortunately, there are
well-known, simple methods for analyzing the required
mapping classes. 
The Euler-Poincar\'e formula (\cite{FLP})
implies that  any \pA\ on the thrice-punctured disk 
has invariant foliations whose only singularities are
one-prongs at the punctures and a one-prong on the outer
boundary. As in the proof of Theorem~\ref{lowerboundthm},
this implies that for any \posp\ $\chi\in\cP_2$,
$h(\chi) = \log(\rho(\Psi(\chi))$ with $\Psi(\chi)$ defined
in \S\ref{anotherrepsec}.

For $j = 1,2$, let $N_j = \Psi(\alpha_j)$.
Now $\chi$ is a word in the generators 
$\alpha_1^{\pm1}, \alpha_2^{\pm1}$, and thus with
exactly the same indices, $\Psi(\chi)$ is
a product of matrices from the set 
$\cN :=\{N_1^{\pm1}, N_2^{\pm1}\}$. 
 Finding the maximal efficiency in $\cP_2$ is then exactly
the same as finding the \gsr\ of $\cN$ (see Remark~\ref{gsree}).
From \eqref{repeq}, 
\begin{equation*}
N_1^{\pm 1} = 
\left(\begin{matrix}
1 & 2 \cdot (-1)^{\pm 1}  \\ 0 &1
\end{matrix}\right), \ \ \
N_2^{\pm 1} = 
\left(\begin{matrix}
1 & 0 \\ 2 \cdot (-1)^{\mp 1} & 1
\end{matrix}\right).
\end{equation*}

The following  proposition is no doubt well known;
we include the simple proof.
\begin{proposition}\label{gsrsimp}
If $\cM =  \{M_1, \dots, M_n\}$ and the matrix $M_1$ has the
property that 
$\|M_1\|_2 = \max\{\|M_{i}\|_2: 1 \leq i \leq n\}$ and further,
$M_1^T \in \cM$, then the generalized spectral radius of
the collection $\cM$ is  $\rho(\cM) = \rho_2(\cM) = 
\|M_1\|_2 = (\rho(M_1 M^T_1))^{1/2}$.
\end{proposition}

\proof\ By \eqref{maxeq}, $\rho(\cM) \leq \|M\|_2$, and since 
$\rho(M_1 M_1^T)^{1/2} = \|M\|_2$, the upper bound is realized
by $M_1 M_1^T$. 
\QED\

\begin{theorem}\label{casetwothm}
The maximal entropy efficiency for \posp\ with
two obstacles is realized
by the protocol $\alpha_1\alpha_2\I$, and so has value
\begin{equation*}
\Eff(2) = \frac{\log(3 + 2\sqrt{2})}{2} = \log(1 + \sqrt{2}). 
\end{equation*}
\end{theorem}
\proof\ The set $\cN$ is closed under taking transposes and
all its elements have the same two-norm, namely 
$(3 + 2\sqrt{2})^{1/2} = 1 + \sqrt{2}$, and
so the theorem follows from Proposition~\ref{gsrsimp}. 
\QED\  

The same argument obviously also shows that the
maximal efficiency is also realized by 
$\alpha_1\I\alpha_2$ or $\alpha_2\alpha_1\I$
or $\alpha_2\I\alpha_1$. Also notice that 
in terms of the usual generators of the braid group
on three strings, $\alpha_1\alpha_2\I = \sigma_1^2  \sigma_1 \sigma_2^{-2}
\sigma_1\I$. Conjugating by $\sigma_1$ we get  
$\sigma_1^2  \sigma_2^{-2}$, and so the path of the
stirrer is a figure eight with an obstacle inside
each loop. This protocol is  called the ``taffy puller''
because it is used in making certain kinds of  candy
(\cf\ Section II in \cite{finnth1}).

\section{The upper bound}\label{upperboundsec}
The main tool for finding an upper bound to 
the entropy efficiency is
the incidence matrix  
of the action of a \posp\ on the fundamental
group of the punctured disk.
 Proposition~\ref{gsrupper}  says that
the \gsr\ of the collection of incidence matrices
gives an upper bound for the entropy efficiency.
 The computation of the appropriate \gsr\ is
possible here because
of the particularly simple form of 
 the action of the generating automorphisms
 $\alpha_k$ on the generating loops $\delta_i$ 
of the free group $\pi_1(D_{N+1}, Q)$ given in \eqref{alphaact}.

\subsection{Incidence matrices and the generalized spectral radius}
We start with a general result that
says that the \gsr\ of the
collection of  incidence matrices is  an upper bound 
for the efficiency of a collection of 
 automorphisms of a free group. 

Assume now that the group  $G = F_m$, the 
free group on the generators $x_1, \dots, x_m$.
Let $\bn_i(g)$ be the number
 of occurrences of the generator $x_i$ or its inverse $x_i\I$ in the
reduced form of $g\in G$. We denote the row vector $\bn(g) = 
(\bn_1(g), \dots, \bn_m(g))$, and so 
$\|\bn(g)\|_1 = \ell(g)$, the word length of $g$.
The \de{incidence matrix} of the automorphism $a:F_m\raw F_m$ 
is the $(m\times m)$-matrix $A$ with 
$A_{ij} = \bn_j(a(x_i))$. 
\begin{proposition}\label{gsrupper}
If $\{a_1, \dots, a_k\}$ is a collection of generators 
of the semigroup $\cA\subset \Aut(F_n)$ whose incidence matrices
are $\{A_1, \dots, A_k\}$, then the entropy efficiency of $\cA$ is
bounded above by the generalized spectral radius of the collection
of incidence matrices or 
\begin{equation*}
\Eff(\cA;\{a_1, \dots, a_n\}) \leq \rho(\{A_1, \dots, A_n\}).
\end{equation*}
\end{proposition}

\proof\ 
Let $a\in\cA$ have incidence matrix $A$.
By construction,
$ \bn(a(g))\leq  \bn(g) A$, and
iterating we have for $n>0$, 
\begin{equation*}
\bn(a^n(g)) \leq \bn(g) A^n .
\end{equation*}
Taking one-norms, $n^{th}$ roots, and then using Gelfand's  formula
\begin{equation*}
\limsup_{n\raw\infty}
\| \bn(a^n(g))\|_1^{1/n} \leq 
\limsup_{n\raw\infty}
\|\bn(g)\|_1^{1/n} \|A^n\|_1^{1/n}
= \rho(A). 
\end{equation*}
Since this is true for all $g\in\F_m$, using the definition
in \eqref{expgrowth} and taking the logarithm 
 we have $h(a) \leq \log(\rho(A))$.
Now  treat $a$ as a product of the generators in reduced form,
$a = a_{i_1} \dots a_{i_k}$. Since incidence matrices ignore
cancellations, $A \geq A_{i_1}\dots A_{i_k}$, and thus because
the matrices are positive,
$\rho(A) \geq \rho(A_{i_1}\dots A_{i_k})$. Therefore
for all $a$ of word length $k$,
\begin{equation*}
\frac{h(a)}{k} \leq \log(\rho( A_{i_1}\dots A_{i_k})^{1/k}),
\end{equation*}
and the result follows.
\QED\

\subsection{The incidence matrices and the form of the upper bound}
 The action of $\alpha_k$ on the generating loops $\delta_i$ 
of the free group $\pi_1(D_{N+1}, Q)$
 as given in \eqref{alphaact} has $(N+1)\times (N+1)$-dimensional
incidence matrix $\bAk = I + S^{(k)}$ with  $S^{(k)}$ defined by 
$S^{(k)}_{kk} = 0$, 
$S^{(k)}_{kj} = 2$ for  $j \not= k$, and  
$S^{(k)}_{ij} = 0$ otherwise.
As in \S\ref{anotherrepsec},  for all $k$, the last row of 
$\bAk$ is $(0\, 0\, \dots 0\, 1)$. Thus if we define
the \de{reduced incidence matrix}  
$\Ak$ as the matrix obtained by deleting the last row
and column of $\bAk$, then 
\begin{equation}\label{chopok}
\rho(\{A^{(1)}, \dots, A^{(N)}\})
= \rho(\{\overline{A}^{(1)}, \dots, \overline{A}^{(N)}\}).
\end{equation}
Thus we may restrict attention to the $\Ak$.  
We still have $\Ak = I + S^{(k)}$, but now
$S^{(k)}$ is $N\times N$. We could have alternatively
defined $\Ak$ as the absolute value of the matrix
$\Psi(\alpha_k) = E^{(k)}$  giving the action
of $\alpha_k$ on lifted chains
to the two-fold cover from \S\ref{anotherrepsec}.

It is easy to check that $\bAk$ is also the incidence
matrix for the functional inverse, $\alpha_k\I$ of $\alpha_k$
and further that $\|\Ak\|_1 = 3$ for all $k$. Thus
by   Proposition~\ref{gsrupper}, \eqref{maxeq} and \eqref{chopok} 
 we have
\begin{proposition}\label{logthree}
If $ \{A^{(1)}, \dots, A^{(N)}\}$ is the set of 
reduced incidence matrices  of the standard set of generators
$\{\alpha_1\pmo, \dots, \alpha_N\pmo\}$ of the group of \posp s
$\cP_N$, then the maximal entropy efficiency in the presence of 
  $N$ obstacles satisfies 
\begin{equation*}
\Eff(N) \leq \rho(\{A^{(1)}, \dots, A^{(N)}\}) \leq \log(3),
\end{equation*}
where $\rho(\{A^{(1)}, \dots, A^{(N)}\})$ is the generalized
spectral radius of the set of incidence matrices.
\end{proposition}

\subsection{Computing the generalized spectral radius}
We can improve the upper bound of $\log(3)$ in
Proposition~\ref{logthree} by computing the
generalized spectral radius 
$\rho(\{A^{(1)}, \dots, A^{(N)}\})$, at least in the sense
of showing that it is realized by an explicit matrix.
\begin{proposition}\label{gsrcompute}
If  $\hHN := A^{(1)} \dots A^{(N)}$, then 
$\rho(\hHN)^{1/N} = \rho(\{A^{(1)}, \dots, A^{(N)}\})$.
\end{proposition}
We shall see in Remark~\ref{notunique} that
the product realizing the \gsr\ in Proposition~\ref{gsrcompute}
is not unique.  

The computation of the 
\gsr\ makes frequent use of the fact that
multiplication by an incidence matrix
$\Ak$ changes the column sum in a fairly
simple fashion. Recall that  $c(M)$ denotes
the column sums of $M$. We then have
\begin{equation}\label{Akacts}
c_j(M\Ak) = c_j(M) + 2 c_k(M) \ \text{for}\ j\not = k,
\ \text{and}\  c_k(M\Ak) = c_k(M).
\end{equation}
 Now the maximal column sum
of a positive matrix is its one-norm, $\| M\|_1 = 
\max c(M)$. Thus if we can determine a sequence of
multiplications by incidence matrices which 
maximizes the column sums, then we know the sequence
which maximizes the one-norm.    
An argument using Gelfand's formula then completes the proof.

Finding the proper sequence is not difficult, again because of
the simple action of the $\Ak$ on column sums. An example
will make it clear. 
\begin{example}
We begin with the identity matrix
$I$ and $c(I) = (1, 1, \dots, 1)$. Since all the entries
are equal, it doesn't matter which $\Ak$ we choose to multiply
by,
and so for simplicity we choose $A^{(1)}$, and  
we get $c(I A^{(1)} ) = (1, 3, 3, 3, \dots, 3)$. Now to
maximize the column sums in the next product, \eqref{Akacts} indicates
we should choose an  index $k$ for which the corresponding
entry in $c(I A^{(1)} )$ is maximal. Again, we may choose
any $k>1$, but for simplicity we choose  $A^{(2)}$ yielding
$c(I A^{(1)} A^{(2)} ) = (4, 3, 6, 6, \dots 6)$. Continuing
we see that at the $j^{th}$-step, the $j^{th}$ entry
of the current column sum will be among the maximal values
and so we choose to multiply by  $A^{(j)}$. It is more or less clear
then that after $N$ steps we will get a maximal column sum
by using the sequence prescribed in $\hHN := A^{(1)} \dots A^{(N)}$ 
\end{example}

We need a few definitions and  preliminary lemmas
to formalize the insight of the example.
Fix $N\geq 3$ and let $\order_N$ be all $N$-tuples of 
 positive integers given in non-increasing order, so 
\begin{equation*}
\order_N= 
\{\vec{a}\in (\Z^+)^n : a_1 \geq a_2 \geq \dots \geq a_N >0 \}.
\end{equation*}
The \de{sorting function} is the map
$S: (\Z^+)^n \to \order_n$ which 
sends a row vector of positive integers to the 
row vector of its entries sorted in non-increasing order.
In keeping with our conventions we will have $S$ 
act from the right, and so $\vb\, S$ is the sorted 
version of $\vb\in (\Z^+)^n$.

It is easy to check that for each $j$ and each $\va\in (\Z^{+})^N$,
$\va\, A_j \in (\Z^{+})^N$, and so we may define 
$W_j: \order_N \raw \order_N$ by
$ \va\, W_j :=  \va\, A_j S$. Using the fact that all entries of
$\va$ are strictly positive we get the formulas 
\begin{equation}\label{wact}
\begin{split}
\va\, W_1 &= (a_2 + 2 a_1, \dots,  a_n + 2 a_1, a_1)\\
\va\, W_j &= (a_1 + 2 a_j, \dots, a_{j-1} + 2 a_j,
a_{j+1} + 2 a_j, a_{n} + 2 a_j, a_j) \ \text{when}\ j>1
\end{split}
\end{equation}
The next lemma collects some useful facts about the 
action of the $W_j$ and the usual partial order on
$(\Z^{+})^N$. Its proof is an elementary exercise using
\eqref{wact}.
\begin{lemma} \label{lem:rel}
Assume $\vec{a},  \vec{b}\in (\Z^{+})^N$.
	
\begin{compactenum}
\item If $\vec{a} \ge \vec{b}$, then 
$\vec{a}\,W_j \ge \vec{b}\,W_j$ for all  $j$.  

\item If $\vec{a} > \vec{b}$, then
$\vec{a}\, W_j >\vec{b}\,  W_j$ for all $j$.      

\item If $j<k$ and $a_j=a_k$, then  
$\vec{a}\, W_j = \vec{b}\, W_k$.

\item If $j<k$  and $a_j>a_k$, then
$\vec{a}\, W_j > \vec{b}\, W_k$.

\item If $j \le k$  and  $\vec{a}>\vec{b}$, then  
$\vec{a}\, W_j >\vec{b}\,  W_k$.
\end{compactenum}
\end{lemma}

The next step is to formalize the various sequences of
multiplying the matrices in our family. 
A \de{strategy} is a list $\us=(s_1, s_2, \dots, s_k)$ 
with $1 \le s_i\le N$ for all $i$. Each strategy 
$\us$ defines a map $W_{\us} : \order_N \to \order_N$ given
by $\va\, W_{\us}:=\va\, W_{s_1} W_{s_2}\dots W_{s_k}$. 

The next definition gives a canonical form
for a strategy which will facilitate comparing the
results of two strategies.
Given $\vec{a} \in \order_N$ 
and a strategy $\us$, the \de{standardization} of $\us$
with respect to $\va$ is the strategy 
$\us'=(s_1', s_2', \dots, s_k')$ defined inductively
as follows.
Let $s_1' = \min\{ i : a_i = a_{s_1}\}$. For
$m>1$, let $b^{(m)} = \va\, W_{s_{1}'} \dots  W_{s_{m-1}'}$, 
and let $s_m' = \min\{ i : b^{(m)}_i = b^{(m)}_{s_m}\}$. 
Note that it follows easily that for all $m$, 
$b^{(m)} = \va\, W_{s_{1}}   \dots   W_{s_{m-1}}$, and
so in particular, 
\begin{equation}\label{standok}
\va\, W_{\us} = \va\, W_{\us'}.
\end{equation} 
If the strategy $\us$ is equal to its standardization
with respect to $\va$, then $\us$ is said to be 
\de{standard} for $\va$. In light of 
\eqref{standok} we henceforth restrict consideration
to standard strategies.

The simplest strategy is to repeatedly use the
first place, which means using the largest entry; for each 
$k$, let $\uonek$ be the $k$-tuple of all ones, 
$\uonek = (1, 1, \dots, 1)$.  The next proposition collects
some useful facts about strategies and products
of the matrices from $\cM$. The first part
says that $\uonek$ it is always the best strategy 
for maximizing the end result.
We adapt the simplified 
notation that given an $k$-tuple $\uell = (\ell_1 \dots \ell_k)$,
$A_\uell = A^{(\ell_1)} \dots A^{(\ell_k)}$.
\begin{lemma}\label{orderstuff}
Assume  $\vec{a} \in \order_n$. 
\begin{compactenum}
\item If $\us = (s_1, \dots, s_k)$ is standard strategy 
for $\va$ and $\us \not = \uonek$, then 
$\va\, W_{\uonek} > \va\, W_{\us}.$
\item If $\us = (s_1, \dots, s_k)$ is any strategy,  
 then $\va\, W_{\uonek} \geq \va\, W_{\us}$.
\item Given any $k$-tuple $\uell = (\ell_1, \dots, \ell_k)$, 
 there exits a strategy 
$\us = (s_1, \dots, s_k)$ 
with $\va\, A_{\uell}\, S = \va\, W_{\us}$.
\item If $\uell$ is the
$k$-tuple $\uell = (1, 2,  \dots k)$ and  
 $\va = (1, 1, \dots, 1)\in (\Z^+)^N$, then
 $\va\, A_{\uell}\, S = \va\, W_{\uonek}$.
\end{compactenum}
\end{lemma}

\proof\
Since by hypothesis, $\us\not= \uonek$, there is
a smallest index  $i$ with $s_i>1$. If $i=1$, let
$\vb = \va$  and if $i>1$, 
let $\vec{b}=\va\, W_{(s_1, \dots, s_{i-1})} = 
W_{\uonekludge}(\va)$.
Since $\us$ is standard for $\vec{a}$, when $i=1$,
$a_{s_i} < a_{1}$ and so  $b_{s_i} < b_{1}$, and
when $i>1$, again since $\us$ is standard for $\vec{a}$,
$b_{s_i} < b_{s_i-1} \leq b_1$. Thus in either
case by Lemma~\ref{lem:rel}(d),
$\vb\, W_1>\vb\, W_{s_i}$ which is to say that 
$\va\, W_{\uonekk}> \va\, W_{s_1, \dots, s_i}.$	
Since $s_n \ge 1$ for all $n>i$, by Lemma~\ref{lem:rel}(e)
$\va\, W_{\uonek} > \va\, W_{\us}$, as required for (a). 

For (b), if $\uonek$ is the standardization of $\us$,
then $\va\, W_{\uonek} = \va\, W_{\us}$,
otherwise use (a).
We prove (c) by induction on $k$, with $k=1$
being trivial. The inductive step is based on the easily
checked observation
that for any (perhaps not ordered vector) $\vd\in (\Z^+)^N$,
and for any $j$, if we let $j'$ be the index with 
$(\vd S)_{j'} = d_j$, then $\vd A_j\, S = \vd\, S A_{j'}\, S$.

Now assume that the result is true for a $k$-tuple $\ell$.
Given $\ell^{+} = (\ell_1, \dots, \ell_k, \ell_{k+1})$,
using the observation of the previous paragraph 
with $\vd = \va\, A_{\ell}$, there is 
a $j'$ with  
$\va\, A_\ell\, A_{\ell_{k+1}}\, S = \va\, A_{\ell}\, S\, A_{j'}\, S 
= \va\, W_{\us}\, A_{j'}\, S$, using the inductive hypothesis.
Thus $\va\, A_{\ell^{+}}\, S = \va\, W_{\us^{+}}$ with
$\us^{+} = (s_1, \dots, s_k, j')$, completing the induction.

The proof of (d) is a simple calculation.
\QED\  

\medskip
\noindent\textbf{Proof of Proposition~\ref{gsrcompute}:}
Given a $k>0$, write $k = m N + j$ with $0\le j <N$ 
and let $\junkk = (\hHN)^m A^{(1)} \dots A^{(j)}$. As a first
step we show
that for any product  of $k$ matrices $P := A^{(j_1)} \dots A^{(j_k)}$,
  $\| \junkk \|_1 \geq \| P\|_1$. 

Let $\va = (1, 1, \dots, 1)\in (\Z^+)^N$.
By Lemma~\ref{orderstuff}(c),  there is a strategy $\us$
 which we may assume to
be standard so that $\va\, P\, S = \va\, W_{\us}$.
Further, by Lemma~\ref{orderstuff}(d),  $\va\, \junkk S$ is produced 
by the strategy $\uonek$, 
so $\va\, \junkk S = \va\, W_{\uonek}$. 
Thus by Lemma~\ref{orderstuff}(b), $\va\, \junkk\, S \geq \va\, P\, S$.
 But for  any positive matrix $M$, $ \max (\va\, M) = \|M\|_1$ and 
  so $\|\junkk\|_1 \geq \|P\|_1$, finishing the proof of
the first step. 

Since $\junkk = (\hHN)^m A^{(1)} \dots A^{(j)}$,
the first step implies that $\|P\|_1 \leq c \|(\hHN)^m\|_1  $, where 
\begin{equation*}
c =  \max\{\|A^{(1)} \dots A^{(j)}\|_1 : 
1\leq j < N \}.
\end{equation*}
But for any matrix $M$, $\rho(M) \leq \|M\|_1$,
and so  for each $k = m N + j$ letting 
$\cM = \{A^{(1)}, \dots, A^{(N)}\}$ we have
\begin{equation*}
\rho_k(\cM) \leq \|(\hHN)^m\|_1^{1/k} c^{1/k} \leq 
\|(\hHN)^m\|_1^{1/(mN)} c^{1/k},
 \end{equation*}
and so using Gelfand's formula,
$\rho(\cM) \leq  \rho(\hHN)^{1/N}$.  
On the other hand, for any $m>0$  
$\rho((\hHN)^m) = \rho(\hHN)^m$, and so
$\rho_{mN}(\cM) \geq \rho(\hHN)^{1/N}$ and thus
$\rho(\cM) \geq \rho(\hHN)^{1/N}$, finishing the proof.
\QED\

\begin{remark}\label{notunique}
The product of matrices in Proposition~\ref{gsrcompute}
is by no means unique.
If $\uell$ is any $N$-tuple which is a cyclic
permutation of $(1, 2,  \dots N)$ and  
 $\va = (1, 1, \dots, 1)\in (\Z^+)^N$, then it is
easy to check that
 $\va\, A_{\uell}\, S = \va\, W_{\uonek}$.
 Thus the joint spectral radius of 
the collection $\{A^{(1)}, \dots, A^{(N)}\}$ it is also
realized by any product of $N$ of the $\Ak$ in which each 
$\Ak$ occurs exactly once.
\end{remark}
	
\subsection{Properties of $\hHN$ and the upper bound}
 Recall that $c_j(M)$ denotes the sum of the entries
in the $j^{th}$ column of the matrix $M$. 
\begin{proposition}\label{hspprop}
The product of incidence matrices $\hHN := A^{(1)} \dots A^{(N)}$
satisfies $c_j(\hHN) = 3^N-2\cdot 3^{j-1} $ 
and $\trace(\hHN) = 3^N -N -1 $, and
so 
\begin{equation*}
\frac{3^N -N -1}{N} \leq \rho(\hHN) \leq 3^N -2.
\end{equation*}
\end{proposition}

\proof\
Fix $N$ and for $1\leq k \leq N$, let 
$\hH^{(k)} = A^{(1)} A^{(2)} \dots A^{(k)}$.
An easy induction yields that $c_j(\hH^{(k)} ) = 3^k - 2\cdot 3^{j-1}$
for $j\leq k$ and $c_j(\hH^{(k)} ) = 3^k$ for $j>k$. The column
sum formula then follows by letting $k=N$. 

For any matrix $M$, it follows from the definition of $A^{(j)}$ that
 $\trace(M A^{(j)}) = \trace(M) + 2 c_j(M) - 2 M_{jj}$.
Since for any $k$, it is easy to check that
 $\hH^{(k-1)}_{k, k} = 1$, and we have
just seen that $c_k(\hH^{(k-1)}) = 3^{k-1}$,   so  another
simple induction   gives $\trace(\hH^{(k)}) = N + 3^k - 2k -1$,
and we again let $k=N$. 

The estimates on the spectral radius then follow from the
general fact that for any $(N\times N)$-matrix $M$,
 $|\trace(M)| /N \leq \rho(M) \leq \|M\|_1$.
\QED\ 

As a consequence of Propositions~\ref{logthree}, \ref{gsrcompute},
and \ref{hspprop} we have our upper bound for the
entropy efficiency.
\begin{theorem}\label{upperboundthm}
The maximal efficiency $\Eff(N)$ of a \posp\  
in the presence of $N$ obstacles satisfies
\begin{equation*}
 \Eff(N) \leq \frac{\log(\rho(\hHN))}{N} \leq \frac{\log(3^N -2)}{N}.
\end{equation*}
\end{theorem}

\bibliography{hsp}
\end{document}